\documentclass[runningheads,a4paper]{llncs}

\usepackage{url}
\urldef{\mailsa}\path|{g.kouyialis14,  r.misener}@imperial.ac.uk|    
\newcommand{\keywords}[1]{\par\addvspace\baselineskip
\noindent\keywordname\enspace\ignorespaces#1}
\usepackage{setspace} 
\usepackage{blkarray}
\usepackage{multirow}
\usepackage{latexsym, times, verbatim, dcolumn, setspace, lscape}
\usepackage{amsmath, amssymb, amsfonts, mathptmx}

\usepackage{amsthm}
\usepackage{textcomp}
\usepackage{verbatim}
\usepackage{upgreek}
\usepackage{epsfig}
\usepackage{float}
\usepackage{bbold}
\usepackage{tabularx}
\usepackage{algorithmicx}
\usepackage[noend]{algpseudocode}
\usepackage[ruled]{algorithm}
\usepackage[T1]{fontenc}
\usepackage[english]{babel}

\usepackage{amsfonts}
\usepackage{color}
\usepackage{tcolorbox}
\usepackage{epstopdf}
\usepackage{xcolor}
\usepackage[margin=25mm]{geometry}
\usepackage{pdfpages}
 \usepackage{sidecap}
\usepackage{wrapfig}
\usepackage{subfig}
\usepackage{tikz}
\usetikzlibrary{matrix,calc}

\usepackage[latin1]{inputenc}
\usepackage[english]{babel}
\usepackage[numbers]{natbib}

\newcommand{\X}{{\bf X}} 
\newcommand{\p}{{\bf p}} 
\newcommand{\Y}{{\mathit{Y}}} 
\newcommand{\E}{{\mathit{E}}} 
\newcommand{\V}{{\mathit{V}}} 
\newcommand{\G}{{\bf G}}

\newcommand{\W}{{\mathit W}}
\newcommand{\Z}{{\mathit Z}}
\newcommand{\A}{{\bf A}}
\newcommand{\B}{{\bf B}}
\newcommand{\Q}{{\bf Q}} 
\newcommand{\x}{{\bf  x}}

\newcommand{\R}{{\mathbb{R}}}
\newcommand\overmat[2]{%
  \makebox[0pt][l]{$\smash{\color{white}\overbrace{\phantom{%
    \begin{matrix}#2\end{matrix}}}^{\text{\color{black}#1}}}$}#2}

\newcommand\partialphantom{\vphantom{\frac{\partial e_{P,M}}{\partial w_{1,1}}}}

\spnewtheorem{exmp}{Example}{\itshape}{\rmfamily}

\spnewtheorem{defn}{Definition}{\bfseries}{\itshape}
\spnewtheorem{prob}{Problem}{\bfseries}{\itshape}
\spnewtheorem{prop}{Proposition}{\bfseries}{\itshape}
\spnewtheorem{method}{Method}{\bfseries}{\itshape}
\spnewtheorem{rem}{Remark}{\bfseries}{\itshape}
\spnewtheorem{graph}{Graph}{\bfseries}{\itshape}
\spnewtheorem{state}{Statement}{\bfseries}{\itshape}

\usepackage{ragged2e}
\usepackage{caption}

\usepackage{multicol, multirow, booktabs, hhline, rotating, array, supertabular}

\def\addsymbol #1: #2#3{$#1$ \> \parbox{5in}{#2 \dotfill \pageref{#3}}\\}

\newlength\savedwidth

\begin{document}

\mainmatter  

\title{Symmetry Detection for Quadratically Constrained Quadratic Programs Using Binary Layered Graphs}

\author{Georgia Kouyialis $\cdot$ Ruth Misener}
\authorrunning{Symmetry Detection for Quadratically Constrained Quadratic Programs Ussing Binary Layered Graphs}

\institute{Department of Computing, Imperial College London, SW7 2AZ, UK\\
\mailsa\\}

\maketitle

\begin{abstract}
Symmetry in mathematical programming may lead to a multiplicity of solutions. In nonconvex optimisation, it can negatively affect the performance of the branch-and-bound algorithm. Symmetry may induce large search trees with multiple equivalent solutions, i.e.\ with the same optimal value. Dealing with symmetry requires detecting and classifying it first. This work develops methods for detecting groups of symmetry in the formulation of quadratically constrained quadratic optimisation problems via adjacency matrices. Using graph theory, we transform these matrices into \textit{Binary Layered Graphs (BLG)} and enter them into the software package \texttt{nauty} \citep{nauty}. \texttt{Nauty} generates important symmetric properties of the original problem.
\end{abstract}

\section{Introduction}\label{sec:intro}
Several geometry problems are mathematically formulated as quadratically constrained quadratic programs \citep{ kucherenko:2007, kallrath:2009}. The occurrence of symmetry in these problems results in many equivalent feasible and optimal solutions that each can be technically generated from the other. Identifying and classifying problem symmetries is an important step towards exploiting tree-based algorithms such as branch-and-cut. This subsequently allows state-of-the-art solver softwares to omit symmetric solutions. This section studies quadratically constrained quadratic programs and the McCormick relaxation method applied on the bilinear terms presented in this formulation. We also provide basic preliminaries on group theory.  

A general formulation of a \textit{Quadratically Constrained Quadratic Program} is given:

\begin{defn} \label{def:QP}
\begin{equation}\label{eq:QP}\tag{$QCQP$}
\begin{aligned} 
\min_{\x\in \R ^n} \quad & f_0(\x) \\
\text{s.t.} \quad & f_k(\x)\leq 0  & \forall \, k=1,\ldots ,m\\
			& x_{i} \in [x_i^{L}, x_{i}^{U}] & \forall \, i = 1,\ldots,n
\end{aligned}
\end{equation}
where 
\begin{equation}
f_k(\x) = \mathop{\sum_{i=1}^{n}\sum_{j=1}^{n}x_i\alpha_{ij}^kx_j} + \sum_{i=1}^{n}\alpha_{i0}^kx_i + \alpha_{00}^k \, \forall \, k=0.\ldots, m
\end{equation}
with coefficients $\alpha^k_{ij} \in \mathbb{R}$ for $i\in \{0,\ldots, n\}$, $j=\{0,\ldots, n\}$ and $k \in\{0,1, \ldots, m\}$ for $x_{i} \in [x_i^{L}, x_{i}^{U}], \, i \in \{1,\ldots, n\}$.
\end{defn}

The property of nonconvexity in mathematical programs imposes more difficulties when we try to solve them. Nonconvexity causes the existence of multiple local optima when we may be seeking a global solution that gives the best optimal value.
There are several convex relaxation techniques for global optimisation problems \citep{Liberti:2004}.
\citet{mccormick:1976} achieves a convex relaxation of quadratically constrained quadratic problems of such problems by adding inequality constraints generated on new auxiliary variables which combine the given ones.
More precisely, a Reformulation Linearisation Technique (RLT) is the McCormick convex and concave relaxation for bilinear terms.

\begin{figure}[h]
\begin{minipage}{0.3\textwidth}\label{fig:RLT}
\small
\includegraphics[width=2.5cm]{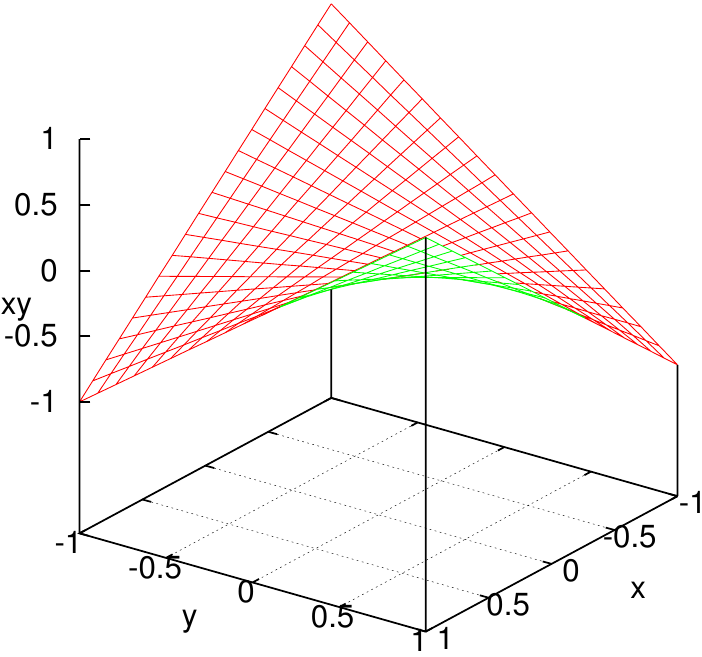}{\\ Bilinear Surface}
\end{minipage}
\begin{minipage}{0.3\textwidth}
\small
\includegraphics[width=2.5cm]{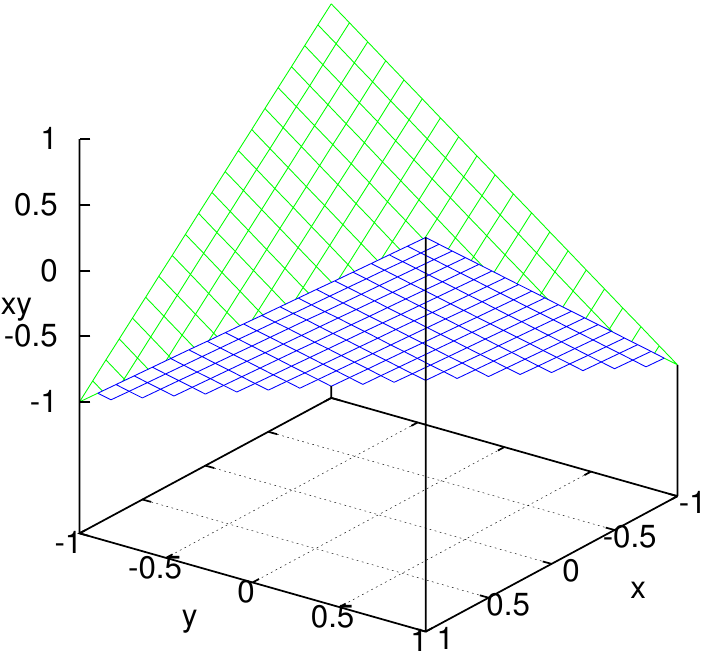}{\\ Lower convex envelope}
\end{minipage}
\begin{minipage}{0.4\textwidth}
\small
\includegraphics[width=2.5cm]{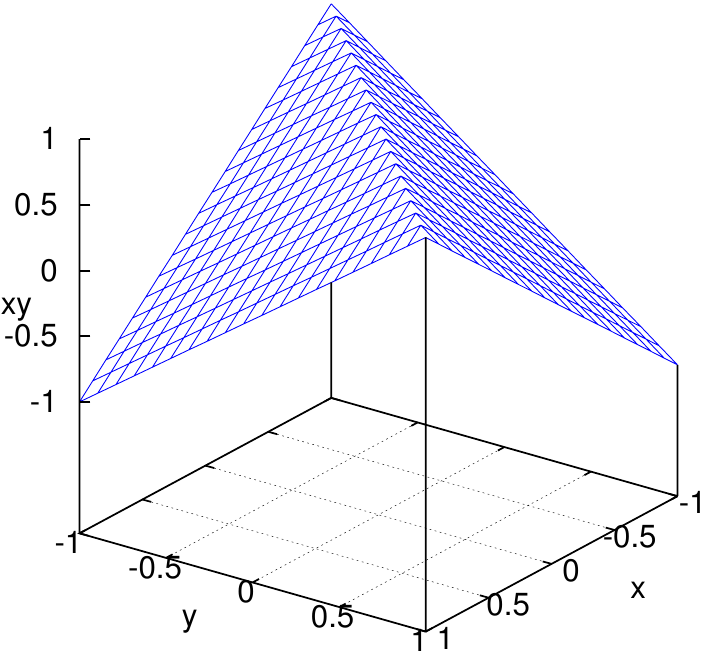}{\\ Upper concave envelope}
\end{minipage}\caption{McCormick convex and concave relaxation for bilinear terms}
\end{figure}

This part follows \citet{anstreicher:2009} and \cite{qualizzaetal-2012} to derive convex relaxation of the original \ref{eq:QP}. For each bilinear term set $X_{ij} = x_i x_j$, the McCormick hull forms under and overestimator constraints, related to the variable bounds $x_i^L \leq x_i \leq x_i^U$, \, $x_j^L \leq x_j \leq x_j^U$
 \[    \left\{ 
  \begin{array}{l}
    X_{ij} \geq x_i^Lx_j + x_j^Lx_i - x_i^Lx_j^L\\
    X_{ij} \geq x_i^Ux_j + x_j^Ux_i - x_i^Ux_j^U\\
      X_{ij} \leq x_i^Ux_j + x_j^Lx_i - x_i^Ux_j^L\\
    X_{ij} \leq x_i^Lx_j + x_j^Ux_i - x_i^Lx_j^U\\
 \end{array} \right.\]
\begin{state} Any quadratic program (QCQP) can be linearised by using the reformulation linearisation technique.
\end{state}
To derive the \citet{mccormick:1976} convex and concave relaxation for bilinear terms, consider any quadratic equation of the form $f_k(\x) = \x^T \Q_k \x + {\p_k}^{T}\x + r_k \leq 0 \,  \forall \, k=\{0,\ldots ,m\}$ and define:

\[{\X} = \x\x^{T} = \left(\begin{array}{c}
x_1\\ \cdot\\ \cdot\\ \cdot\\ x_n \end{array} \right) \left(\begin{array}{ccccc}
x_1 & \cdot & \cdot &\cdot & x_n \end{array}\right) = \left(\begin{array}{ccccc}
x_1x_1 & x_1x_2 & \cdot &\cdot & x_nx_n\\
x_2x_1 & \cdot & \cdot & \cdot &\cdot \\
\cdot & \cdot &\cdot &\cdot &\cdot\\
\cdot & \cdot &\cdot &\cdot &\cdot\\
x_nx_1 & \cdot &\cdot &\cdot & x_nx_n \end{array}\right)\]

Rewrite each quadratic expression using the inner product:

\begin{equation*} \label{ref}
\x^{T} \Q_k \x = \mathop{\sum_{i=1}^{n}\sum_{j=1}^{n}\Q_{kij}x_ix_j} = \Q_k \bullet {\X} = \mathop{\sum_{i=1}^{n}\sum_{j=1}^n \Q_{kij}X_{ij}}
\end{equation*}

Now use the variable bounds of $x_i$,$x_j$ to obtain the constraints of the following linearised optimisation problem \ref{eq:LQ}.
Note that for $i\neq j$, $X_{ij} = X_{ji}$ as the matrices are symmetric above their diagonal, and $\X = {\X}^T$. Hence we can define the linear form of \ref{eq:QP} as

\begin{defn}\label{def:LQ}
\begin{equation} \label{eq:LQ}\tag{$LQP$}
\begin{aligned} 
\max \quad &  \Q_0\bullet {\X} + {\p_0}^{T}\x + {r_0} \\
\text{s.t.} \quad &  \Q_k \bullet {\X} + {\p_k}^{T}\x + r_k \leq 0 & \forall \, k=\{1,\ldots ,m\}\\
            & {\X} - \x^L\x^T - \x(\x^L)^T \geq -\x^L(\x^L)^T\\
            & {\X} - \x^U\x^T - \x(\x^U)^T \geq -\x^U(\x^U)^T\\
            & {\X} - \x^L\x^T - \x(\x^U)^T \leq -\x^L(\x^U)^T\\
            & {\X} - \x^U\x^T - \x(\x^L)^T \leq -\x^U(\x^L)^T\\
            & {\X} = {\X}^T \\
            & \x\in[\x^{L},\x^{U}] 
\end{aligned}
\end{equation}

where $\x\in \mathbb{R}^n$, $\Q_0, \ldots, \Q_m \in \mathbb{R}^{n\times n}$, are $n$ by $n$ matrices and $\p_k \in \mathbb{R}^n$, $r_k \in \mathbb{R}$.
\end{defn}

Group theory studies the algebraic structure of a wide range of objects with or without special properties. An important class of groups is the one of transformations (symmetries). We use the class of permutation groups, e.g.\ the symmetric and the cyclic group, to describe the transformations in geometric objects \citep{armstrong:1988, cameron:1999}.

A \textit{group} $(\W, \cdot)$ is a nonempty set $\W$ with a binary operation $\cdot$ on $\W$ satisfying the following properties: If $g, \, z$ $\in \W$, then $g\cdot z$ is also in $\W$; $g\cdot (z\cdot d)=(g\cdot z)\cdot d$ for all $g, \, z, \, d \in \W$; Every group $\W$ has an identity element $I$, such that $g\cdot I = I\cdot g = g, \; for each g \in \W$; If  $g \in \W, \, \exists$ $g^{-1} \in \W$ such that $g\cdot g^{-1} = g^{-1}\cdot g = I$.
A \textit{subgroup} $\Z$ of a group $\W$ is a nonempty subset of $\W$ that forms a group itself under the operation induced by $\W$.
Two groups $\W$, $\Z$ are \textit{isomorphic} if $\exists$ a bijective function $\phi: \W \rightarrow \Z$ that satisfies: $\phi(I) = I$ for an identity element $I$; $\phi({g^{-1}}) = {\phi }(g)^{-1},\ \forall g \in \W$; $\phi(gz) = \phi(g)\phi(z), \ \forall g,\ z \in \W$. A \textit{group automorphism} is an isomorphism from a group to itself. 
A \textit{permutation} of a set $\textit{Y} = \{1, \, \ldots, \, n\}$ is a bijective function ${\bf \pi: \textit{Y} \longrightarrow \textit{Y}}$. 
A \textit{permutation group} ${\Pi}^n$  is a finite group whose elements are permutations of a given set \textit{Y} and whose group operation is composition of permutations in the group.
For permutations $\pi \in {\Pi}^n$, $\sigma \in {\Pi}^m$, $\A(\pi,\sigma)$ is a matrix obtained by permuting the columns of $\A$ by $\pi$ and the rows of $\A$ by $\sigma$.
The \textit{symmetric group} $S^n$ is the group of all permutations of a given set \textit{Y}.
In mathematics, representation theory studies ways to represent the elements of groups as linear transformations of vector spaces.
The \textit{symmetry group} of an object is the group of all transformations under which the object is invariant with group operation the composition of such transformations. 
A \textit{cyclic group} is a group that can be generated by a single element e.g.\ $\textbf{$C_n$} = \langle r | r^n = 1, n \in \Z \rangle$ is the cyclic group of order $n$.
An object has cyclic symmetries if it is invariant under the transformations in a cyclic group. 

A main part of this paper evolves around the symmetry in the original nonconvex \ref{eq:QP} problem and how it is affected after relaxing the problem which leads to an ordinary linear program \ref{eq:LQ}.
\begin{table}\caption{Table of Notation.}\label{tbl:notation}
\centering
\begin{tabular}{ l  l  l  l}
\toprule
 Symbol & Description & Symbol & Description\\
\midrule
$x_i$ & Variables  & $I$ & Identity element\\
$\x$ & Vectors of variables  & $  \pi, \sigma$ & Permutations\\
$\alpha$ & Coefficient  & $\Pi^n$ & Set of all permutations\\
$ \bf c, \bf b, \bf p$ & Vectors of parameters  & $S_n $ & Symmetric group order n\\
$\A, \Q$ &  Matrices of parameters  & $\Y$ & Sets\\
$ \X$ & Matrix of auxiliary variables & $f, h, \phi$ & Functions\\
$ M, IM, JM, KM$ & Sparse representations of matrices & $\G, \bf H$ & Graphs\\
$\cal{F}$ & Set of feasible solutions & $\E, \V$ & Set of edges, vertices\\
$\cal{G}$, $\tilde{{\cal{G}}}$ & Symmetry groups & $e$ & Edges in the graph\\
$\W, \Z$ & Groups & $u, v$ & Nodes in the graph\\

\bottomrule
\end{tabular}
\end{table}
This paper proceeds as follows: 
Section \ref{sec:background} surveys the relevant literature and provides a specific geometry problem which motivates the study of symmetry.
Section \ref{sec:def_symmetry_QCQP} formally defines symmetry in quadratically constrained quadratic optimisation problems and identifies the role of integrality and nonconvexity in such cases. Section \ref{sec:data_structures_detect_symmetry_DAG} evolves around the formulation symmetries in optimisation problems and the graph structures that currently exist in literature for detecting such symmetries. 
Section \ref{sec:data_structures_matrices} suggests two different methods on forming a problem as an adjacency matrix and explains how to convert these matrices into graphs.
Section \ref{sec:data_structures_BLG} introduces \textit{binary layered graphs} and describes how to use the proposed matrices and construct these graphs. Section \ref{sec:data_structures_computational_case} shows how to automatically detect symmetry using software package \texttt{nauty}. This work concludes in Section \ref{sec:data_structures_conclusion} with a discussion on the proposed structures and comparison to other methods.
\section{Symmetry in Mathematical Programming}
\label{sec:background}
\subsection{Literature Review}
\citet{margot:2010} defines symmetry in Integer Linear Programming (ILP) of the form\\ $P^{L} = \min_{\x\in \mathbb Z ^{n}}\{{\bf c}^{T}\x \mid \A \x \geq \bf b \}$ where $\A \in \R ^{m \times n}$ and vectors $\bf {c}$ $\in \R^ n$, $\bf b$ $\in \R ^{m}$. Symmetry is the set of variable permutations under which any feasible solution remains feasible and the objective function value is invariant. 
Let ${\cal{F}}$ be the set of feasible solutions of any problem $P^{L}$.
The symmetry group is: 
\[\tilde{{\cal{G}}}(P^L) = \{\pi \in {\Pi}^n \mid \forall \quad \hat{\x} \in {\cal{F}}, \quad \pi (\hat{\x})\in {\cal{F}}  \quad \text{and} \quad {\bf c}^T\pi (\hat{\x}) ={\bf c} ^T\hat{\x}\}\]
\citet{liberti:2012} studies and extends the definition of symmetry to mixed-integer nonlinear optimisation problems.
 
Symmetric structure in optimisation may be viewed through the lens of group theory for \ref{eq:QP} \citep{bodi:2013, herr:2013}.
 In many situations though, it is difficult to detect the symmetry of the original problem and its polyhedron representation \citep{bremner:2014}. 
 The formulation group is a subgroup of the symmetry group and reflects the symmetric properties of the variables and the constraints of an optimisation problem. 
 Symmetry handling approaches present methodologies to associate optimisation problems with graph representations from which the graph automorphism is generated by using software tools \citep{puget:2005, berthold:2008, margot:2010, libertiref:2012, bodi:2013, knueven:2017}.
 
Many researchers exploit the above information and the insights of several problems; covering problems \citep{margot:2003, margot:2007}, scheduling and packing problems \citep{anjos:2010, costa:2013} and engineering problems as the unit commitment problem and heat exchanger network synthesis \citep{ostrowski-etal:2015, magnago:2016, kouyialis:2016}. 
They identify the presence of symmetry and in some cases propose symmetry handling approaches for problems with known symmetric structure. 
The improved performance of the solvers validates the efficiency of these techniques \citep{rehn:2015}.
 However, they are problem specific and cannot be generalised to other problems.
 
There are several methods to exploit symmetry which are categorised as static and dynamic methods.
Static methods adjoin new constraints to the formulation in order to make some symmetric optima infeasible. 
Sherali and co-workers add symmetry breaking constraints or perturb the objective function \citep{sherali:2001, ghoniem:2011}.
Other researchers investigate the orbitopes of a problem \citep{berthold:2008, kaibel:2009, kaibel:2011}: convex hull of 0-1 matrices that represent possible solutions to packing and partitioning constraints. 
The new constraints yield to a reformulation which is guaranteed to keep at least one symmetric optimum feasible.
Orbitopes have additionally been considered for cutting planes \citep{friedman:2007, hojny:2015}.
\citet{liberti:2008} automatically generates symmetry handling inequalities, whereas other works study inequalities which exploit multiple variable orbits \citep{ostrowski:2014, dias:2015}; the groups of variables that can be sent to each other under some actions (permutations in the group) which  are equivalent with respect to symmetry of the problem.

In the dynamic category are approaches which modify the solution method i.e.\ the search tree algorithm to recognise and exploit symmetry dynamically as it goes along. For example, constraints can be derived for each node in the tree to forbid the isomorphic nodes \citep{gent:2000, gent:2005, ramani:2005}.
Another way of exploiting symmetry in B$\&$B is given by isomorphism pruning \citep{ margot:2002, margot:2003a, margot:2003} and orbital/constrained orbital branching \citep{ostrowski:2008, ostrowski:2009}.
By introducing artificial variables, \citet{fischetti:2017} reformulate the problem to a reduced problem which considers only variables of symmetry orbits instead of all variables, so-called orbital shrinking.

While most of these works consider the symmetry representation as a step enclosed by the scope of handling symmetry, \citet{liberti:2012} is the first who stated the importance of a practical and general representation of symmetry. 
He uses expression trees to explicitly capture the structure of an optimisation problem and develops the ROSE \citep{rose} reformulation software engine that produces a file representation of the problem as \emph{Directed Acyclic Graphs} (\emph{DAG}).

The work introduced in this thesis concerns with the improvements on symmetry detection, which is the first phase of symmetry handling techniques. 
Symmetry representation is an elementary process given to the software package \texttt{nauty}, on which all the following steps to break symmetry depend.
 Hence it is very essential to guarantee and increase its correctness and efficiency.

 \subsection{Motivation}
To isolate this phenomenon we present a prototype circle packing problem. Visually consider a problem of locating two identical circles $(c_1, c_2)$ with centre coordinates $(x, y)$, $(x', y')$ in a unit square. Figure 2 illustrates this problem.
The optimisation problem is to make the circles as large as possible without overlapping. There are four ways to locate these circles and they are related by rotations and reflections. Mathematically speaking, there are four sets of feasible (approximated) solutions which give the same objective value; distance between their centre coordinates \citep{pack}.
\begin{figure}[t!]
\centering
\begin{tikzpicture}[every node/.style={draw,circle,minimum size=0.57cm,inner sep=0pt}]
\draw (0,0) -- (1,0) -- (1,1) -- (0,1) -- cycle;
\node[draw=none] at (0.2, 0.8){\scriptsize $\bf 1.$};
\node[ fill=blue!5] at (0.293, 0.293) {};
\node[fill=yellow!5] at (0.707, 0.707) {};
\draw[->] (1.1,0.5) to (1.9, 0.5) ;
\node[draw=none] at (1.5, 0.7) {\scriptsize rotation};

\draw[->] (1.1,-0.3) to (1.9, -1) ;
\node[draw=none] at (1.5, -0.1) {\scriptsize reflections};
\draw[->] (1.9,-0.3) to (1.1, -1) ;

\node[draw=none] at (2.8, 0.8){\scriptsize $\bf 2.$};
\draw (2,0) -- (3,0) -- (3,1) -- (2,1) -- cycle;
\node[ fill=blue!5] at (2.293, 0.707) {};
\node[fill=yellow!5] at (2.707, 0.293) {};
\draw[->] (3,-0.1) to (3, -0.9) ;
\node[draw=none] at (2.5, -0.5) {\scriptsize rotation};

\node[draw=none] at (0.8, -1.2){\scriptsize $\bf 4.$};
\draw (0,-1) -- (1,-1) -- (1,-2) -- (0,-2) -- cycle;
\node[fill=yellow!5] at (0.293, -1.293) {};
\node[ fill=blue!5] at (0.707, -1.707) {};
\draw[->] (0,-0.9) to (0, -0.1) ;
\node[draw=none] at (0.5, -0.5) {\scriptsize rotation};

\node[draw=none] at (2.2, -1.2){\scriptsize $\bf 3.$};
\draw (2,-1) -- (3,-1) -- (3,-2) -- (2,-2) -- cycle;
\node[fill=yellow!5] at (2.293, -1.707) {};
\node[ fill=blue!5] at (2.707, -1.293) {};
\draw[->] (1.9,-1.5) to (1.1, -1.5) ;
\node[draw=none] at (1.5, -1.3) {\scriptsize rotation};
\end{tikzpicture}\caption{Example that shows four different ways of locating two circles in a unit square which lead to the same optimal minimum distance between their centre.}
\end{figure}
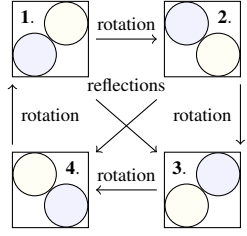\label{Fig:circ}

 As shown in Figure 2: 
 \begin{equation*}
 \begin{split}
  1 = \{(0.293, 0.293), (0.707 , 0.707)\}\\ 2= \{(0.293, 0.707), (0.707, 0.293)\}\\ 3= \{(0.707, 0.707), (0.293, 0.293)\}\\ 4= \{(0.707, 0.293), (0.293, 0.707)\} \end{split}
  \end{equation*}  However only one can be considered as unique as all the others can be obtained by permuting the variables of the problem.
Consider solution $1$ and permute variables $(y \, y')$ to get solution $2 = \{(0.293, 0.707), (0.707, 0.293)\}$, $(x \, x')$ to get $4 = \{(0.707, 0.293), (0.293, 0.707)\}$, apply both permutations $(y \, y')(x \, x')$ to get solution $3 = \{(0.707, 0.707), (0.293, 0.293)\}$. Permutations $(x \, y)$ and/or $(x' \, y')$ take solution $1$ to itself \citep{costa:2013}.
The exchange of the variables of the problem which leaves the set of feasible solutions and the objective function value unaffected is the \emph{symmetry} in an optimisation problems.
In a branch-and-bound framework, symmetry is an optimal solution with different configuration  leading to the same objective function value.
 In worst case, $B\&B$ exhaustively enumerates all feasible solutions. 
 Hence, the presence of symmetry can cause unexpectedly large trees which immediately affects the time that is taken for the algorithm to terminate and the problem to be solved.
Hence exploiting symmetry is a challenge. 
Identifying and classifying problem symmetries is an important step towards exploiting tree-based algorithms such as branch-and-cut. 
This subsequently allows state-of-the-art solver software to omit symmetric solutions. 

\section{Symmetry Group of Quadratically Constrained Quadratic Programs}\label{sec:def_symmetry_QCQP}

After surveying the available sources for detecting symmetry, we contemplate the explanation of symmetry given by \citet{margot:2010} which is also presented by \citet{liberti:2012} on different problems. 
Under a set of permutations of the problem variables, each feasible solution can be mapped to another solution having the same value and the whole set of feasible solutions ${\cal{F}}$ can be mapped to itself.

Modifying this definition to the case of quadratic problems we define symmetry in QCQP:

\begin{defn}{Symmetry group}\\
$$\tilde{{\cal{G}}}(P^Q) = \{\pi \in {\Pi}^n \mid \forall \quad \hat{\x} \in {\cal{F}}, \quad  \pi(\hat{\x})\in {\cal{F}}  \quad \text{and} \quad f_0(\pi(\hat{\x})) = f_0( \hat{\x}) \}$$

\end{defn}

The symmetry group is based on the feasible set of solutions of an optimisation problem. 
Deriving this set is impractical in our work.
Hence, the scope of this paper is to efficiently associate data structures with optimisation problems which can generate the formulation group: a set of permutations that fix the problem formulation. \citet{liberti:2012} proves that the formulation group is a subset of the symmetry group. Definitions and relevant structures of this group are provided in the next section.

This work discusses the symmetry in nonlinear QCQPs, as well as the symmetry in linearised cases of \ref{eq:LQ}, which arise after applying the RLT to the QCQP formulation.

Relaxing the problem using the RLT technique adds constraints which may alter the symmetric structure of the problem. 
We evaluate and state how the integrality and nonlinearities affect the symmetry group of the original problem. The following examples show that another challenging aspect of detecting symmetry is the fact that symmetry in the original problem does not imply the same symmetry in the relaxed problem and vice versa.

 
 \begin{exmp}\label{ex:symmetry_QCQP_ex1}
 \begin{align*} 
\min \quad & x_1 + x_2   \\
\text{s.t.} \quad & 2{x_1} + x_2 \leq 2\\
  & x_1, x_2 \in \{0,1\}
\end{align*}
 \end{exmp}
The set of feasible solutions is $\{(0,0), (1,0), (0,1)\}$. Hence, the symmetry group is characterised by permutation $\pi = (x_1, x_2)$ under which any feasible solution remains within the set of feasible solutions and the objective function value is invariant. On the other hand, if we relax the integrality constraints over a continuous range $x,y \in [0,1]$ the feasible solution $(0.5,1)$ under permutation $\pi(0.5,1) = (1,0.5)$ violates the linear constraint.  
Hence, if the integrality restrictions in a Mixed-integer QCQP are relaxed and the symmetry group is defined over the feasible set of solutions, then it is not necessarily a subgroup of the symmetry group in the relaxation.
The next example shows that the symmetry group of the relaxation is also not a subgroup of the symmetry group in the original problem.

 \begin{exmp}\label{ex:symmetry_QCQP_ex2b}
	 \begin{align*}  
\min \quad & x_1 + x_2    \\
\text{s.t.} \quad &  x_1 + x_2 \leq 1\\
  &  x_1 \in [0,1]\\
&  x_2 \in \{0,1\}
\end{align*} 
\end{exmp}

Similar to Example \ref{ex:symmetry_QCQP_ex1} permutation $\pi = (x_1, x_2)$ is the symmetry of the relaxed problem, but solution $(0.5,1)$ is not feasible under this permutation for the original problem. 

Another example shows that the original optimisation problem might not inherit any symmetry. But under relaxation, e.g. \ McCormick, symmetries arise.

\begin{exmp}
	\begin{equation}\label{ex:sym1}\tag{$QP$}
\begin{aligned} 
\min \quad & -x_1 - x_2^2 \\
\text{s.t.} \quad
            & x_2 + x_3 \geq 1\\
            &  x_1 \geq 2x_3 -1\\
            & x_1 \leq x_3\\
            & x_1, x_2, x_3 \in \{0,1\}            
\end{aligned} 
\end{equation}


	 	\begin{equation}\label{ex:sym2}\tag{$LQP$}
\begin{aligned} 
\min \quad & -x_1 - x_4 \\
\text{s.t.} \quad
             & x_2 + x_3 \geq 1\\
            &  x_1 \geq 2x_3 -1\\
            & x_1 \leq x_3\\                    
            &  x_4 \geq 2x_2 -1\\
            & x_4 \leq x_2 \\            
             & x_1, x_2, x_3, x_4 \in \{0,1\}   
\end{aligned}
\end{equation}

\end{exmp}

For the symmetry group the feasible set of solutions:\\  ${\cal{F}}({QP}) = \{ (1\, 1\, 1), (0\, 1\, 0), (1\, 0\, 1)\}$ with ${\tilde{\cal{G}}}({QP}) = \{ I\}$.\\ 
${\cal{F}}(LQP) = \{(1\, 1\, 1\, 1), (1\, 0\, 1\, 0), (0\, 1\, 0\, 1)\}$ with\\ ${\tilde{\cal{G}}}({LQP}) = \{I, (x_1\, x_3), (x_2\, x_4), (x_1\, x_3)(x_2\, x_4), (x_1\, x_2)(x_3\, x_4), (x_1\, x_4)(x_2\, x_3)\}$ hence, ${\cal{G}}(LPQ) \nleq {\cal{G}}(QP)$.

\section{Formulation Symmetry Detection via Directed Acyclic Graphs}
\label{sec:data_structures_detect_symmetry_DAG}

The formulation group of a mathematical optimisation problem is defined by \citet{libertiref:2012} as the set of permutations of the variable indices for which the objective function and the constraints are the same. Hence, for \ref{eq:QP}:

\begin{defn}{Formulation group of \ref{eq:QP}} 
$${{\cal{G}}}(P^Q) = \{\pi \in {\Pi}^n \mid \quad f_0(\pi(\x )) = f_0(\x) \quad \text{and} \quad \exists \quad \sigma \in {\Pi}^m (\sigma f_k(\pi (\x )) = f_k(\x)) \}$$
\end{defn}

\subsection{Graph Theory}\label{sub:graph_theory}

A \textit{graph} is a tuple $\G = (\V, \E)$ where $\V$ is a (finite non-empty) set of vertices and $v \in \V$ is called a vertex and $\E \subset \V \times \V$ is a finite collection of edges and $e = \{u,v\} \in \E$ is called an edge.
An edge from a vertex to itself $e = \{u\}$ is said to be a \textit{loop}.
A \textit{weighted graph} $\bf K$ is a triplet $\bf K = (\V, \E, w)$ where $w : \E(\bf K) \rightarrow \mathbb{R}$.
A $\ell$-colouring of a labelled graph $\G = (\V, \E, c)$ is a function $c: V(\G) \rightarrow \{0,1,\ldots, \ell-1\}$ where $k$ is the number of colours. The vertices of one colour form a colour class and $\G$ is defined as a \textit{vertex coloured graph}. 
Two simple graphs $\G = \{\V(\G), \E(\G)\}$ and $\bf H = \{\V(\bf H), \E(\bf H)\}$, are \textit{isomorphic graphs}, denoted $\G \cong \bf H$, if
$\exists$ a bijective function $f : \V(\G) \rightarrow \V({\bf H})$, such that for each edge $\{u\, , v\} \in \E(\G)$ there is an edge $\{f(u), f(v)\} \in \, \E(\bf H)$. Under this relation, any set of adjacent vertices $\E(\V) = \{\{u, v\} | u, v \in \V, u \neq v\}$ remains adjacent.  
An automorphism is an isomorphism of a graph to itself.
Given a graph $\G$, a permutation $\pi$ of $\V (\G)$ is a \textit{graph automorphism} of $\G$, 
if $\forall \, u, \, v \, \in \, \V (\G)$ there exist an edge $\{u\, , v\} \in \E(\G)$ then under any permutation $\pi$ remains in the set of edges as $\{\pi(u),  \pi(v)\} \in \, \E(\G)$.
Consider pairs $(\G,c)$ where $\G$ is a graph and $c: V(G) \rightarrow \{0,\dots,\ell - 1\}$ is a $\ell$-colouring of $\G$. A \textit{colour - preserving isomorphism} from $(\G, c)$ to $(H, c')$ is a bijection $\pi : V(\G) \rightarrow V(\bf H)$ such that $\pi$ is an isomorphism from $\G$ to $\bf H$ and $c(v) = c'(\pi(v))$ $\forall v \in V(\G)$. 
A colouring of the vertices is also referred to as a partition, and the colour classes as the cells of the partition.
A \textit{graph partitioning} of $\G$ into $\ell$ parts is a collection of nonempty disjoint subsets $\V_0, \ldots , \V_{\ell-1}$ for $\ell \in \mathbb{Z}$ whose union is $\V$, i.e.\, $\V = \V_0 \cup \ldots  \cup \V_{\ell-1}$ $\forall \ell$.

\subsection{Expression Trees}\label{sub:data_structures_exptrees}
To compare two functions, \citet{libertiref:2012} suggests to compare their expression trees. An expression tree as first introduced by \citet{Crawford:1996} and explained by \cite{ramani:2005} for Constrained Programming (CP) is used to represent algebraic functions, since it can visually present the structural relation of its components.  
To guarantee that a tree correctly represents an algebraic expression, it should contains all of the component i.e.\ operations, constants and variables. Therefore, tree nodes are categorised into three types: operator nodes, constant nodes, and variable nodes. All the actions to modify an expression tree, like removing parentheses and merging similar terms, are in accordance with the laws of algebra. The rank of a node $v$ is the maximum number of edges taken to reach a node and all the leaf nodes are of rank zero. 
The basic rules are:
\begin{enumerate}
\item Operators are distinguished in: binary (difference, power) and $k$-ary (sum and product) for positive integer $k$.
\item Leaf nodes: labelled with variable symbols and numerical constants.
\item Non-leaf nodes: labelled with operator symbols.
\end{enumerate}

\begin{figure}[h]
\centering
    \includegraphics[width=0.5\textwidth]{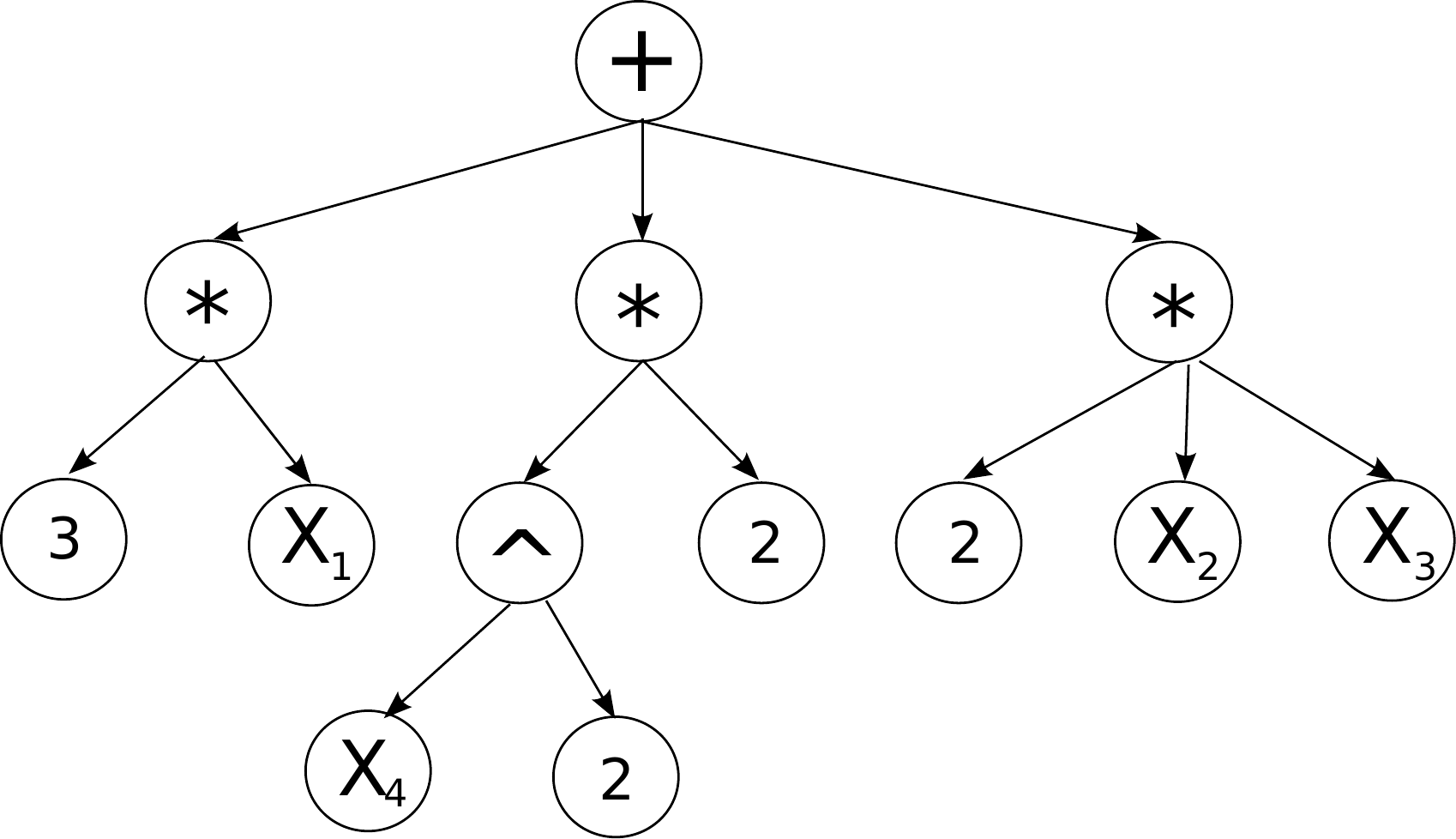}\caption{Example of expression tree representation.}\label{fig:exptree}
  \end{figure}

A simple example is provided in Figure \ref{fig:exptree}:
\begin{exmp}
$$3 x_1+2 x_4^2+2 x_2 x_3$$
\end{exmp}

An algorithm for the tree comparison starts at the root node by comparing their attributes and values. If the current nodes are equal, then it descends down to the child nodes and carries on the comparison in the same way. The comparison steps are recursively executed until the test function reports the existence of equivalence or detects that the two trees violate an equivalence criteria.

Designing equivalence test functions seems reasonable, however, it is not practical. 
Such tests might require a large number of numerical comparisons, and so they would be algorithmically intractable. 
To validate the correctness of this method,  \citet{libertiref:2012} claims that  $f_1, f_2 : \mathbb{R}^n \rightarrow \mathbb{R}$, are equivalent if they have the same range of feasible domain, i.e.\ $dom(f_1) = dom(f_2)$ and $\forall x \in dom(f_1) f_1(x) = f_2(x)$.
Based on which he proves that the formulation group of a mathematical program is a subset of its symmetry group, i.e.\ ${\cal{G}(P)} \leq \tilde{{\cal{G}}}(P)$.

Moreover, in terms of the problem formulation, the role of convex relaxation is highly significant. The definition must be strictly applied in every nonlinear term otherwise the symmetric properties of the problem can be affected as in the following example.
\begin{exmp}\label{ex:symmetry_QCQP_ex4} 
Consider the original problem:
\begin{align*} \label{ex:form}
\min \quad & -x_1^2 - x_2^2 \\
\text{s.t.} \quad & 0 \leq x_1 \leq 1 \\
            &  0 \leq x_2 \leq 1
\end{align*}
\end{exmp}
with $\tilde{{\cal{G}}}(P^Q) = {\cal{G}}(P^Q) = \{I, (x_1, x_2)\}$\\

If we generate RLT constraints only for $x_1^2$, let $x_4 = x_1x_1$ then we get $(P^{LQ})$ which has no symmetric properties,  
$\tilde{{\cal{G}}}(P^{LQ}) = \{I\}$. 
 Following the definition of convex relaxation,\\ let $x_4 = x_1x_1$,  $x_3 = x_2x_2$ then we get $P^{LQ'}$ with ${{\cal{G}}} = \{I, (x_1, x_2)(x_3, x_4)\}$. 

\subsection{Directed Acyclic Graphs}\label{sub:data_structures_DAG}
\citet{liberti:2012} reduces the formulation symmetry problem to the graph isomorphism problem. As an extension of the expression tree structures for single functions, he introduces a coloured \textit{DAG} for multiple functions appearing in mathematical programs. Such functions have the same variable (argument) list, so the trees can share the same variable leaf nodes. Further simplifications for duplicated nodes and algebraic equivalence are applied by \citet{libertiref:2012}.

 A major component of these structures is an equivalence relation on the graph vertices which determines the interchangeability of two vertices. Subsequently, a graph colouring partitions the \textit{DAG} vertices and identifies the subsets of nodes which can be permuted. 

The vertex set of an expression graph is partitioned according to the following rules:
\begin{enumerate}
\item Root nodes that represent the constraints can be permuted iff they have the same RHS.
\item Variable nodes can be permuted, iff they are of the same type and same range.
\item Constant nodes can be permuted, iff they have the same rank level and value. 
\item Operator nodes  can be permuted, iff they have the same rank level and value.
\item The order of a child node can not be exchanged iff the operator node is non-commutative. 
\end{enumerate}

Two important theoretical results support the correctness of \textit{DAG} constructions.
\citet{ramani:2005} prove that:
\begin{theorem}
The symmetries of the constraints of the given mathematical problem, correspond one-to-one to the symmetries of the graph.
\end{theorem}

\citet{libertiref:2012} proves how to map the automorphism group of a \textit{DAG} graph to the formulation group of the original problem. 
\begin{theorem}
A subgroup of the automorphism group of a \textit{DAG} that fixes the variable nodes of the graph is equivalent to the formulation group of the original problem.
\end{theorem}

\section{Symmetry Representation via Matrices}
\label{sec:data_structures_matrices}

This section proposes and discusses structures to detect the formulation group which captures the symmetric nature of a given linear and nonlinear programming problem. 

\subsection{Matrix Structures}\label{sub:data_structures_matrix}
This part suggests two different methods of forming a problem as a matrix. 
The definition of the formulation group of a problem depends on these matrices. 
We transform each matrix into a graph for detecting and classifying the automorphism group which reveals the symmetry of the original problem.
The presence of linear and bilinear terms in quadratic problems though indicates their difficulty. 
Consider the formulation of \ref{eq:QP} with functions 
\begin{equation*}
f_k(\x) = \mathop{\sum_{i=1}^{n}\sum_{j=1}^{n}x_i\alpha_{ij}^kx_j} + \sum_{i=1}^{n}\alpha_{i0}^kx_i + \alpha_{00}^k \, \forall \, k=0.\ldots, m
\end{equation*}
with coefficients $\alpha^k_{ij} \in \mathbb{R}$ for $i\in \{0,\ldots, n\}$, $j=\{0,\ldots, n\}$ and $k \in\{0,1, \ldots, m\}$ for $x_{i} \in [x_i^{L}, x_{i}^{U}], \, i \in \{1,\ldots, n\}$.

\begin{method}\label{meth:2}       
Create a tensor $\A^Q \in \mathbb{R}^{(n+1)\times (n+1)\times (m+1)}$ with entries $a_{ij}^k$ as shown in Figure 4.
\end{method}

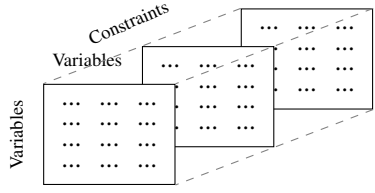
\begin{figure}[h!]
        \centering
\begin{tikzpicture}
\def\xs{1.3} 
\def\ys{0.5} 
\def\nm{3} 
\foreach \x in {1,2,...,\nm}
{

\matrix [draw, 
         fill=white, 
         ampersand replacement=\&] 
(mm\x)
at(-\x * \xs, -\x * \ys) 
{   
 \node { ...}; \& \node { ...}; \& \node { ...};\\
    \node { ...}; \& \node { ...}; \& \node {... };\\
        \node {... }; \& \node { ...}; \& \node { ...};\\
                \node { ...}; \& \node { ...}; \& \node { ...};\\
};
}
\node[draw=none, yshift=2ex, xshift=0em, anchor=west] at (mm3.north west) {\scriptsize Variables};
\node[draw=none, yshift=0ex, xshift=-1em, rotate=90] at (mm3.west) {\scriptsize Variables};
\node[draw=none, yshift=2ex, xshift=-0.5em, rotate=24] at (mm2.north west) {\scriptsize Constraints};
\draw [dashed,gray](mm1.north west) -- (mm\nm.north west);
\draw [dashed,gray](mm1.north east) -- (mm\nm.north east);
\draw [dashed,gray](mm1.south east) -- (mm\nm.south east);
\end{tikzpicture}\caption{Matrix illustration of Method \ref{meth:2}.}
\end{figure}\label{fig;tensor}

Each matrix corresponds to a QCQP equation and the rows and columns to a constant element and the variables of the QCQP, capturing the relations between the bilinear term. The following example shows this idea.
\begin{exmp}
	\begin{equation}\label{ex1d}\tag{$QP_1$}
\begin{aligned} 
\max \quad & 3x_1 +  3x_4 + 2x_2x_3 \hspace{1.9cm} (c_0)\\
 	    & x_2 + {x_1}^2 + 1\leq 0\hspace{2.1cm} (c_1)\\
            & x_3 + x_4^2 + 1\leq 0 \hspace{2.27cm} (c_2)\\
            & x_2 + x_3 + 1\leq 0   \hspace{2.27cm} (c_3)\\
            & x_1, x_2, x_3, x_4 \in [0,1]         
\end{aligned}
\end{equation}

${\A}^{Q_1}$ = 
{\scriptsize
\begin{align*}
A_{c_0} =
 \begin{pmatrix}
  0 & 3 & 0 & 0 & 3\\
  3 & 0 & 0 & 0 & 0\\
  0 & 0 & 0 & 2 & 0\\
  0  & 0  & 2 & 0 & 0 \\
   3  & 0  & 0 & 0 & 0
 \end{pmatrix}
 A_{c_1} =
 \begin{pmatrix}
  1 & 0 & 1 & 0 & 0\\
  0 & 1 & 0 & 0 & 0\\
  1 & 0 & 0 & 0 & 0\\
  0  & 0  & 0 & 0 & 0 \\
  0  & 0  & 0 & 0 & 0  \end{pmatrix}
 A_{c_2} =
  \begin{pmatrix}
  1 & 0 & 0 & 1 & 0\\
  0 & 0 & 0 & 0 & 0\\
  0 & 0 & 0 & 0 & 0\\
  1  & 0  & 0 & 0 & 0 \\
  0  & 0  & 0 & 0 & 1 
 \end{pmatrix}
 A_{c_3} =
 \begin{pmatrix}
  1 & 0 & 1 & 1 & 0\\
  0 & 0 & 0 & 0 & 0\\
  1 & 0 & 0 & 0 & 0\\
  1  & 0  & 0 & 0 & 0 \\
  0  & 0  & 0 & 0 & 0 
 \end{pmatrix}
\end{align*}
}
\end{exmp}

%
%

Consider the problem \ref{eq:LQ} which incorporates the constraints of the original problem and the RLT constraints formed by McCormick relaxation for each nonlinear term. Let $\hat{m} = (1 + m + (\#\text{ of non linear terms})\times 4)$ and $\hat{n} =  (1 + n + \#\text{ of non linear terms})$.
\begin{method}\label{meth:1} 
Create a 2 dimensional matrix  $\A^{LQ} \in \mathbb{R}^{\hat{m} \times \hat{n}}$, with entries the coefficients of \ref{eq:LQ} $a_{kj}$ as shown in Figure 5.
\end{method}
\begin{figure}[!h]
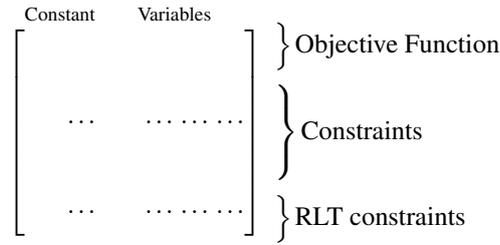

    \centering
        \centering
$\begin{matrix}

 \begin{bmatrix}
 \overmat{Constant}{} \qquad \qquad \, \,
     \overmat{Variables}{} &\\[0.3em]\\
     \ldots  & \ldots & \ldots & \ldots \\[0.5em]
   &  \\[0.5em]
    \ldots & \ldots & \ldots & \ldots  \\[0.5em]
    
  \end{bmatrix}
  \begin{aligned}
  &\left.\begin{matrix}

  \partialphantom  \\[0.5em]
  \end{matrix} \right\} %
  \text{Objective Function}\\
  &\left.\begin{matrix}
\partialphantom  \\[0.5em]
  \partialphantom  \\[0.5em]
  \end{matrix} \right\} %
  \text{Constraints}\\
  &\left.\begin{matrix}
 
\partialphantom  \\[0.5em]

  \end{matrix}\right\}%
  \text{RLT constraints}\\
 \end{aligned}
 \end{matrix}$\caption{Matrix illustration of Method \ref{meth:1}.}   
\end{figure}    

The number of columns set as $\hat{n}$ consists of a constant element, each variable and the auxiliary variables introduced for nonlinear terms. 
The number of rows say $\hat{m}$ consists of the objective function and all the constraints of the problem. 
Note that the maximum number of nonlinear terms is: $\frac{n(n+1)}{2}$. 
Consider the linearised form of $QP_1$ by introducing the auxiliary variables, $X_{23} = x_2x_3$, $X_{11} = x_1^2$, $X_{44} = x_4 ^2$ and add the McCormick relaxation constraints.
Adding the RLT constraints lead to the following example:
\begin{exmp}
	\begin{equation}\label{ex2d}\tag{$LQP_1$}
\begin{aligned} 
\max \quad & 3x_1 + 2X_{23} + 3x_4 \hspace{1.9cm} (c_0)\\
\text{s.t.} \quad & X_{11} + {x_2} + 1 \leq 0 \hspace{2cm} (c_1)\\
            & x_3 + X_{44} + 1\leq 0 \hspace{2cm} (c_2)\\
            & x_2 + x_3 + 1\leq 0 \hspace{2.2cm} (c_3)\\            		
	        & x_2 + x_3 - X_{23} - 1 \leq 0  \hspace{1.37cm} (c_4) \\ 
      		& X_{23} - x_2 \leq 0 \hspace{2.6cm} (c_5) \\
			& X_{23} - x_3 \leq 0 \hspace{2.6cm} (c_6) \\      		
      		& 2x_1 - X_{11} - 1\leq 0  \hspace{1.9cm} (c_7) \\ 
      		& X_{11} - x_1 \leq 0  \hspace{2.6cm} (c_8) \\      		
      		& 2x_4 - X_{44} - 1 \leq 0 \hspace{1.9cm} (c_9) \\ 
      		& X_{44} - x_4 \leq 0 \hspace{2.6cm} (c_{10}) \\
      		& X_{23}, X_{11}, X_{44} \geq 0 \quad \\ 
		& x_1, x_2, x_3, x_4 \in [0,1]  
\end{aligned}
\end{equation}

${\A}^{LQ_1}$ = 
{\scriptsize
\[
\begin{blockarray}{(cccccccc)c}
  I & x_1 & x_2 & x_3 & x_4 & X_{11} & X_{23} & X_{44} & \\
  0 & 3 & 0 & 0 & 3 & 0 & 2 & 0 & C_0 \\
  1 & 0 & 1 & 0 & 0 & 1 & 0 & 0 & C_1 \\
  1 & 0 & 0 & 1 & 0 & 0 & 0 & 1 & C_2 \\
  1 & 0 & 1 & 1 & 0 & 0 & 0 & 0 & C_3 \\
  -1 & 0 & 1 & 1 & 0 & 0 & -1 & 0 & C_4 \\
  0 & 0 & -1 & 0 & 0 & 0 & 0 & 0 & C_5 \\
  0 & 0 & 0 & 0 & 0 & 0 & 1 & 0 & C_6 \\
  -1 & 2 & 0 & 0 & 0 & -1 & 0 & 0 & C_7 \\
  0 & -1 & 0 & 0 & 0 & 1 & 0 & 0 & C_8 \\
  -1 & 0 & 0 & 0 & 2 & 0 & 0 & -1 & C_9 \\
  0 & 0 & 0 & 0 & -1 & 0 & 0 & 1 & C_{10} 
\end{blockarray}\]}
\end{exmp}


If we compare these two methods we observe that Method \ref{meth:1} has potentially considerably fewer entries. Consider the case where a problem has $n$ variables. Method \ref{meth:1} requires 4 new constraint for each nonlinear term. Then in the worst case scenario there are $(1+n)(1+ \frac{n}{2})(1+m+2n^2+2n)$ entries in contrast to the Method \ref{meth:2} which has $(1+n)(1+n)(1+m)$. In most cases, Method \ref{meth:2} has fewer entries than Method \ref{meth:1}. There exist pathological cases, e.g.\ fully dense formulations as $m > 3n^2 + 6n + 3$, where Method \ref{meth:1} has fewer entries. The graph transformation is based on the number of entries of these matrices. Hence, dealing with smaller graphs reduces their complexity and the procedure time taken to generate their symmetric properties and to compare them.


Next, we define the formulation group of a matrix; necessary for detecting symmetry. 

\begin{defn}{Formulation group of the matrix $\A^{LQ}$} 
$${\cal{G}}(\A^{LQ}) = \{\pi \in {\Pi}^{\hat{n}} \mid \exists \sigma \in {\Pi}^{\hat{m}} \quad \text{such that} \quad \A(\sigma, \pi) = \A\}$$
\end{defn}

The set of permutations of the columns of $\A^{LQ}$ such that there is a corresponding permutation of the rows that when applied yields the original matrix.
For permutations $\pi \in {\Pi}^n$, $\sigma \in {\Pi}^m$, $\A(\pi,\sigma)$ is a matrix obtained by permuting the columns of $\A$ by $\pi$ and the rows of $\A$ by $\sigma$.

\begin{defn}{Formulation group of the matrix $\A^{Q}$}
$${\cal{G}}(\A^{Q}) = \{\pi \in {\Pi}^n \mid \exists \sigma \in {\Pi}^m \quad \text{such that} \quad \A(\pi, \pi, \sigma) = \A\}$$
\end{defn}

The set of permutations of the columns and rows of each matrix in $\A^{Q}$ under which the matrix yields to its original form. The same permutation $\pi$ acts both on rows and columns of the matrix $\A^{Q}$ which represent the same number and type of variables. 

Matrix representation shows that exchanging the columns and rows of a matrix which subsequently means changing the position of the variables and constraints of the problem, leads to an equivalent problem. \citet{margot:2002, margot:2003a} proves that the formulation group of similar matrices is a subset of the symmetry group of the original problem, i.e ${\cal{G}}\leq \tilde{{\cal{G}}}$. 
The corresponding set of permutations that fixes the problem formulation is called the problem symmetry. \citet{liberti:2008} proposes an automatic way to compute permutations of the formulation group and proves that is a subset of the symmetries in the solution group of a MILP in the general form. As far as we know, this is the only approach for automatic symmetry detection that does not reduce the problem to a graph. Computational experiments report that finding elements of the symmetry automatically is too costly in terms of CPU time.

\subsection{Converting Matrices to Edge-Labelled Vertex-Coloured Graphs}\label{sub:data_structures_converting}
The matrix representations proposed in this paper include all the elements of an optimisation problem by construction: variables, constraints and coefficients. 
The main idea of this work is to convert such matrices into edge-labelled vertex-coloured graphs associated with the basic elements of the problem and then map the graph automorphisms to the original problem symmetries. \cite{margot:2010} states that mapping the instance of a problem to a coloured graph is a standard procedure \citep{ramani:2005, salvagnin:2005, ramani:2006}. Colour preserving automorphisms of such graphs correspond to problem symmetries. 
A similar idea on how to convert the matrices in Section \ref{sub:data_structures_matrix} to edge-labelled vertex-coloured graphs is given here. 

Since many of the matrix values ${\bf A}$ are 0, a sparse matrix representation of Method \ref{meth:2} is used to reduce space in memory and time accessing all the coefficient of the problem. Consider a tuple $\A = (\bf M, \bf I, \bf J, \bf K)$ of vectors $\bf M, \bf I, \bf J, \bf K$ $\in \R^s$ with maximum size $s = (n+1)(n+1)(m+1)$. 

\begin{figure}[h]
\centering
\includegraphics[scale=1.1]{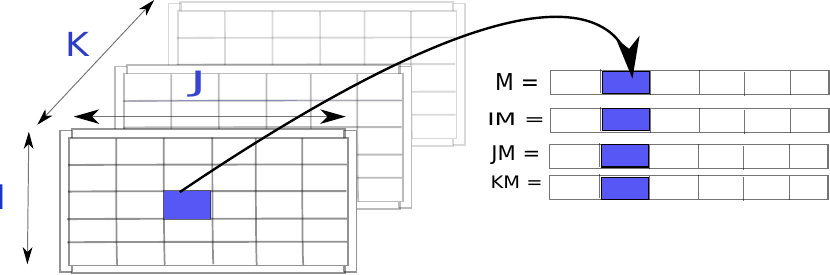}\label{def:sparse}\caption{Sparse matrix representation.}
\end{figure}

\begin{itemize}
\item $\bf {M}$ = $(M_1, \ldots, M_s)$ is a vector with all non zero entries of a matrix stored from left to right and from top to bottom.
\item $\bf J $ = $ (J_1, \ldots, J_s)$ represent the column indices correspond to non zero entries.
\item  $\bf I $ = $ (I_1, \ldots, I_s)$ represent the row indices correspond to non zero entries.
\item $\bf K $ = $ (K_1, \ldots, K_s)$ represent the matrix (constraint) indices correspond to non zero entries. 
\end{itemize}

Recall Example \ref{ex1d}. For Method 1, since each matrix is symmetric, without loss of generality we consider only the upper triangular matrices of each case. 
\begin{itemize}
 \setlength{\itemsep}{0pt}
  \setlength{\parskip}{0pt}
  \setlength{\parsep}{0pt}
  \item $ \bf M_{QP_1} $ = $ (3 \, 3\, 2 \, 1\, 1\, 1\, 1\, 1\, 1\, 1\, 1\, 1)$
\item $     \bf I_{QP_1} $ = $ (0 \, 0 \, 2 \, 0  \, 0 \, 1 \, 0 \, 0 \, 4 \, 0 \, 0 \, 0 \,)$ 
\item $     \bf J_{QP_1} $ = $ (1  \,4 \, 3 \, 0 \, 2  \, 1  \, 0  \,  3 \,  4 \, 0 \, 2 \, 3)$
\item $     \bf K_{QP_1} $ = $ (0 \, 0 \, 0 \, 1\, 1\, 1\, 2\, 2\, 2\, 3\, 3\, 3)$
\end{itemize}
Note that this is not the case for Method 2 since the matrix illustration of Example \ref{ex2d} is not symmetric.

Using these vectors, we construct edge-labelled vertex-coloured graphs which are variants of constraint/variable incidence graphs. 
Consider a graph $\G = (\V, \E, c)$ corresponding to an instance $\bf M, I, J, K$. 
The function ${c}: \E \rightarrow r$, for $r\in \{0,\ldots, \ell -1\}$ is an edge colouring and $\ell \in \mathbb{Z}^+$ is the unique number of different coefficients in $\bf M$. Each of these unique elements in vector $\bf M$ is stored in a vector $\bf {U}$ $\in R^n$.
The vertex set is partitioned (coloured) into four subsets as explained in Section 4.3.2, $\V_F$ a set containing a node for the objective function, $\V_{C}$ nodes for the constraints,  $\V_S$ a constant node and $\V_R$ nodes for the variables. Note that the definition of the automorphism with respect to colours states that each vertex can only be mapped onto a vertex of the same colour. 

\begin{figure}[h]
\centering
\begin{tikzpicture}[every node/.style={draw,circle,minimum size=0.5cm,inner sep=0pt}]

\node[draw=none] (t5) at (0.75,2.35) {};
\node[draw=none] (t6) at (0.75,3.25) {};
\node[draw=none] (t4) at (0.75,2.25) {};
\node[draw=none] (t3) at (0.75,1.35) {};

\node[orange] (S7) at (-1,-0.2)  {\scriptsize $C_{m}$};
\node[draw=none] (S6) at (-1,0.2)  {\scriptsize $\bullet$};
\node[draw=none] (S5) at (-1,0.5)  {\scriptsize $\bullet$};
\node[draw=none] (S4) at (-1,0.8)  {\scriptsize $\bullet$};
\node[orange] (S3) at (-1,1.31)  {\scriptsize $C_{2}$};
\node[orange] (S2) at (-1,2.31)  {\scriptsize $C_{1}$};
\node[red, opacity=0.4] (S1) at (-1,3.31)  {\scriptsize $C_{0}$};

\node[blue] (D7) at (1,-0.2)  {\scriptsize $X_{n}$};
\node[draw=none] (D6) at (1,0.2)  {\scriptsize $\bullet$};
\node[draw=none] (D5) at (1,0.5)  {\scriptsize $\bullet$};
\node[draw=none] (D4) at (1,0.8)  {\scriptsize $\bullet$};
\node[blue] (D3) at (1,1.31)  {\scriptsize $X_{2}$};
\node[blue] (D2) at (1,2.31)  {\scriptsize $X_{1}$};
\node[teal] (D1) at (1,3.31)  {\scriptsize $X_{0}$};

\node[draw=none] (s2) at (0.23,2.31)  {};
\node[draw=none] (s3) at (0.23,1.31)  {};

\node[draw=none] (d1) at (1.25,3.31) {};
\node[draw=none] (d2) at (1.25,2.31)  {};
\node[draw=none] (d3) at (1.25,1.31)  {};

\draw[->, purple] (S3)  edge node[draw=none,rectangle,minimum size=0cm,above]{\scriptsize } (D2);
\draw[->, purple] (S1)  edge node[draw=none,rectangle,minimum size=0cm,above]{\scriptsize } (D2);
\node[draw=none] at (-0.1,3) {\scriptsize \color{purple}{$w_1$}};
\draw[->] (S1)  edge node[draw=none,rectangle,minimum size=0cm,above]{\scriptsize } (D3);
\draw[->, yellow] (S2)  edge node[draw=none,rectangle,minimum size=0cm,above]{\scriptsize } (D1);
\node[draw=none] at (-0.7,2.66) {\scriptsize \color{yellow}{  $w_3$}};
\draw[->] (S2)  edge node[draw=none,rectangle,minimum size=0cm,above]{\scriptsize} (D7);
\node[draw=none] at (0.2,2.31) {\scriptsize \color{black}{  $w_2$}};
\draw[->, brown] (S7)  edge node[draw=none,rectangle,minimum size=0cm,above]{\scriptsize} (D3);
\node[draw=none] at (-0.8,0.2) {\scriptsize \color{brown}{ $w_{\ell}$}};


\end{tikzpicture}\caption{Weighted graph representation for Method \ref{meth:1}.}\label{gr:wgraph1}
\end{figure}
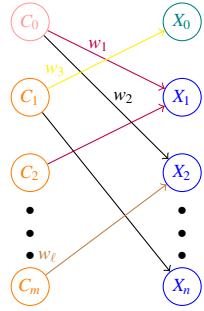

For Method \ref{meth:1}, we construct the following edge coloured graph with edge set initially empty $\E = \emptyset$. For $i = \{0,\ldots, s\}$ where $s = |I| = |K| = |M|$ add an edge  $v_{I_i}^{(r)}$ to $v_{K_i}^{(r)}$, i.e.\ from a vertex in the set that represents the constant element / variables to a vertex in the set of the objective function / constraints, with the relevant colour as shown in Figure \ref{gr:wgraph1}.

For Method \ref{meth:2} the graph construction also accounts edges between nodes in the variable set to show the bilinear relations and loops for quadratic terms. Initially for graph $\G = (\V, \E, c)$ let $\E = \emptyset$. For $i = \{0,\ldots, s\}$:
\begin{itemize}
\item If $I_i = J_i$ then $\E = \E \cup \{\{(v_{I_i}, v_{K_i})^r\} \cap \{(v_{I_i})^r\}\}$ 
\item else for $I_i \neq J_i$ then $
						\E = \E \cup \{\{(v_{I_i}, v_{K_i})^r\}
						 \cap \{(v_{J_i}, v_{K_i})^r\}
							 \cap \{(v_{I_i}, v_{J_i})^r\}\}
					$
\end{itemize}

\begin{figure}[h]
\centering
\begin{tikzpicture}[every node/.style={draw,circle,minimum size=0.5cm,inner sep=0pt}]

\node[draw=none] (t5) at (0.75,2.35) {};
\node[draw=none] (t6) at (0.75,3.25) {};
\node[draw=none] (t4) at (0.75,2.25) {};
\node[draw=none] (t3) at (0.75,1.35) {};

\node[orange] (S7) at (-1,-0.2)  {\scriptsize $C_{m}$};
\node[draw=none] (S6) at (-1,0.2)  {\scriptsize $\bullet$};
\node[draw=none] (S5) at (-1,0.5)  {\scriptsize $\bullet$};
\node[draw=none] (S4) at (-1,0.8)  {\scriptsize $\bullet$};
\node[orange] (S3) at (-1,1.31)  {\scriptsize $C_{2}$};
\node[orange] (S2) at (-1,2.31)  {\scriptsize $C_{1}$};
\node[red, opacity=0.4] (S1) at (-1,3.31)  {\scriptsize $C_{0}$};

\node[blue] (D7) at (1,-0.2)  {\scriptsize $X_{n}$};
\node[draw=none] (D6) at (1,0.2)  {\scriptsize $\bullet$};
\node[draw=none] (D5) at (1,0.5)  {\scriptsize $\bullet$};
\node[draw=none] (D4) at (1,0.8)  {\scriptsize $\bullet$};
\node[blue] (D3) at (1,1.31)  {\scriptsize $X_{2}$};
\node[blue] (D2) at (1,2.31)  {\scriptsize $X_{1}$};
\node[teal] (D1) at (1,3.31)  {\scriptsize $X_{0}$};

\node[draw=none] (s2) at (0.23,2.31)  {};
\node[draw=none] (s3) at (0.23,1.31)  {};

\node[draw=none] (d1) at (1.25,3.31) {};
\node[draw=none] (d2) at (1.25,2.31)  {};
\node[draw=none] (d3) at (1.25,1.31)  {};

\draw[->, purple] (S1)  edge node[draw=none,rectangle,minimum size=0cm,above]{\scriptsize } (D3);
\draw[->, purple] (S1)  edge node[draw=none,rectangle,minimum size=0cm,above]{\scriptsize } (D2);
\node[draw=none] at (-0.1,3) {\scriptsize \color{purple}{$w_1$}};
\draw[->, yellow] (S2)  edge node[draw=none,rectangle,minimum size=0cm,above]{\scriptsize } (D1);
\node[draw=none] at (-0.7,2.66) {\scriptsize \color{yellow}{  $w_2$}};
\draw[->, yellow] (S2)  edge node[draw=none,rectangle,minimum size=0cm,above]{\scriptsize} (D7);
\draw[->, purple] (D2)  edge node[draw=none,rectangle,minimum size=0cm,above]{\scriptsize} (D3);
\draw[->, brown] (S7)  edge node[draw=none,rectangle,minimum size=0cm,above]{\scriptsize} (D3);
\draw[->, brown] (S3)  edge node[draw=none,rectangle,minimum size=0cm,above]{\scriptsize} (D3);
\node[draw=none] at (-0.8,0.2) {\scriptsize \color{brown}{ $w_{\ell}$}};
\draw[->, yellow] (D1)  to [bend left,looseness=1.2] node[draw=none,rectangle,minimum size=0cm,above]{\scriptsize} (D7);
\draw[->, brown] (D3) to [out=10,in=300,looseness=5] (D3);


\end{tikzpicture}
\caption{Weighted graph representation for Method \ref{meth:2}.}\label{gr:wgraph2}
\end{figure}
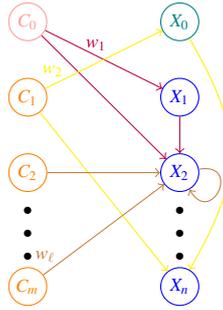


\section{Formulation Symmetry Detection via Binary Layered Graphs}
\label{sec:data_structures_BLG}

To detect symmetry, we use software \texttt{nauty} \citep{nauty}. 
Ideally, each variable could be a vertex in the graph and each coefficient a label of an edge connecting the vertices involved using edge vertex coloured graphs in Section \ref{sub:data_structures_converting}. But \texttt{nauty} \citep{nauty} accepts only vertex coloured graphs, so we associate edge colours with layers in a graph. 
\citet{nauty} state that graphs similar to Section \ref{sub:data_structures_converting} and matrix representations in Section \ref{sub:data_structures_matrix} are isomorphic to a general vertex coloured layered graph representation. According to this statement we describe how to illustrate an optimisation problem with \textit{binary layered graph} (\textit{BLG}) structures (\cite{nauty}). 
In this paper, we convert the adjacency matrices in Section \ref{sub:data_structures_matrix} into \textit{binary layered graphs} (\textit{BLG}) and generate the automorphism group of such graphs that projects the symmetry in the original optimisation problems.

\citet{nauty} explain how to convert a graph $\G = (\V, \E, c)$ with colouring ${c}: \E \rightarrow \{0,\ldots, \ell -1\}$ of $\ell$ colours into an \textit{$\ell$ - layering graph}. 
First, replace  each vertex $v_j \in \V$ with a fixed connected graph of $\ell$ vertices $v_j^{(0)},\dots,v_j^{(\ell -1)}$. If an edge $(v_j,v_{j'})$ has colour $r$, add an edge from $v_{j}^{(r)}$ to $v_{j'}^{(r)}$. Finally, partition the vertices by the superscripts, $V_r = \{ v_0^{(r)},\dots, v_{n-1}^{(r)}\}$.

\subsection{Binary Layered Graph Representation}\label{sec:data_structures_BLG_transformation}

\begin{figure}[h]
\centering
    \includegraphics[width=0.4\textwidth]{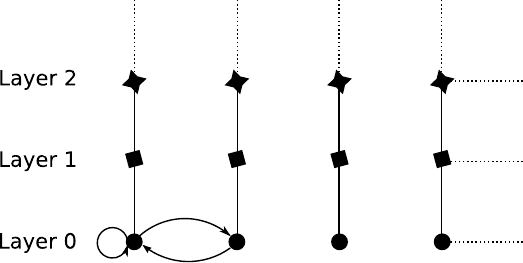}\caption{General form of a Binary Layered Graph.}\label{fig:BLG}
  \end{figure}

We use a binary representation to avoid many layers in $\G$ when the number of colours is large. 
\begin{defn}{Binary Layered Graph}\label{def:BLG}\\
Let $\ell$ be the number of edge labels of $\G$. 
A \textit{BLG} is an edge-labelled vertex-coloured graph $\B$. 
Each vertex colour is associated with a binary representation. 
The number of layers of $\B$ is:
\begin{equation}\label{eq:log}
L = \lceil{\log _2 \left( \ell + 1 \right)} \rceil \quad \text{for} \quad \ell \in \mathbb{Z}
\end{equation}
\end{defn}
Assign a unique positive integer $\mu (z)$ to each  unique element $z$ in vector $\bf {U}$.
The set $\{\mu (z) \, |\, z\in \bf U \}$ is a set of edge labels for $\B$.  
For each $\mu (z)$ compute a binary representation.
\begin{equation}\label{eq:binary}
z = c_{L-1} 2^{L-1} + c_{L-2}2^{L-2} + \ldots + c_{0}2^0, \quad \text{for}\quad c_t \in \{0,1\} \, t = \{0, \ldots, L-1\}
\end{equation}
For nonzero $c_{t}$, the powers of $2$ reveals which layers encode that value.
If $c_t = 1$, add a new edge from $v_i^{(t)}$ to $v_j^{(t)}$ for every $c_t \in \{c_1,\ldots,c_{L-1}\}$.
The form of a layered graph is shown in Figure \ref{fig:BLG}.

In \ref{eq:QP} it is important to consider both nonlinear terms and how to incorporate the different variable coefficients in the graph representation. \textit{BLG} structures can handle this situation with loops as described below. This section shows how to illustrate different mathematical problems as graphs.

Next we describe the different graphic illustrations of problems with a finite number of algebraic expressions. 

\subsection{Graph Structures}\label{sub:graph} 

Section \ref{sub:data_structures_matrix} presents two methods on how to associate matrices with optimisation problems \ref{eq:LQ} and \ref{eq:QP}. 

For $i = \{1, \ldots, n\}$, $k = \{1, \ldots, m\}$, $\ell \in \mathbb{Z}$ where $n$ is the number of variables, $m$ is the number of constraints and $\ell$ the number of unique coefficients in each problem (P). 
The following graph representation skeletons are presented for $\G_{P} = (\V_{P}, \E_{P})$. The vertex set consists of $V_{F}$ set containing vertices associated with the objective function, $V_{C}$ with the constraints and $V_{S}$ with a constant and $V_{R}$ the variables. 
Similar to \cite{liberti:2012}, we define an equivalence relation $\sim$ on $V_{P}$ as follows:
\begin{equation*}
\begin{aligned}
\forall u,v \in V_{P} \, u\sim v \, \implies & (u,v \in V_{F} \wedge \ell (u) = \ell (v))\\
& \vee (u,v \in V_{C} \wedge \ell (u) = \ell (v))\\
& \vee (u,v \in V_{S} \wedge \ell (u) = \ell (v))\\
& \vee (u,v \in V_{R} \wedge \ell (u) = \ell (v))  
\end{aligned}
\end{equation*}
i.e.\ vertices on the same vertex set and layer are in the same partition and can be exchanged.

\begin{graph}\label{gr:1} represents linear problems (originally or after applying RLT) and matrix representation in Method \ref{meth:1}.
\end{graph}
The number of layers $L = \lceil{\log _2 \left( \ell + 1 \right)} \rceil + 1$. The total number of vertices is: $|\V| = (\hat{n}+1)(L - 1) + \hat{m} + 1$.
The vertex set consists of (layer 0) vertices that correspond to the objective function and the constraints of the problem. 
On every other layer, there are copies of these nodes as shown by the vertical lines. Then on the top layer there is one vertex for a constant element and vertices for each QCQP variable.
From nodes in (layer 0) and its copies, we add edges with endpoints the nodes on the top layer, based on which variable is included on each constraint and what is the coefficient in front of this variable. 

\begin{figure}
    \centering
      {\includegraphics[width=0.5\textwidth]{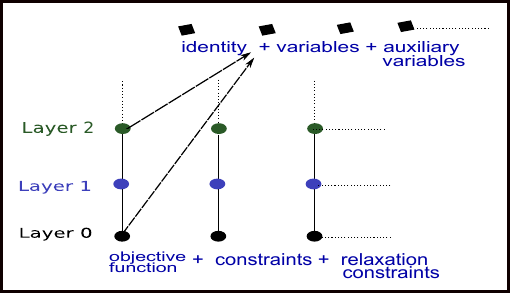}}
    \caption{Illustration of \textbf{Graph 1}.}\label{fig:graph1}
\end{figure}


\begin{graph}\label{gr:3} represents quadratic (nonlinear) $QCQP$ problems and matrix on Method \ref{meth:2}. 
\end{graph}
The graph consists of two different parts with number of layers $L = \lceil{\log _2 \left(\ell + 1 \right)}\rceil + 2$; the vertices for the objective function and each constraint and layers of copies of these constraints (connected with vertical edges). In this part the horizontal edges encode the coefficients of the problem. 
The total number of vertices is $|\V| = 2(n + 1) + (m + 1)(L - 2)$.
On the upper part as shown in Figure \ref{fig:BLG}, there are vertices for a constant element and each variable and a layer of copies of variables (connected with vertical edges). On this layer the horizontal edges and loops distinguish the relations of linear and bilinear terms.
\begin{figure}[h]
    \centering
   \includegraphics[width=0.5\linewidth]{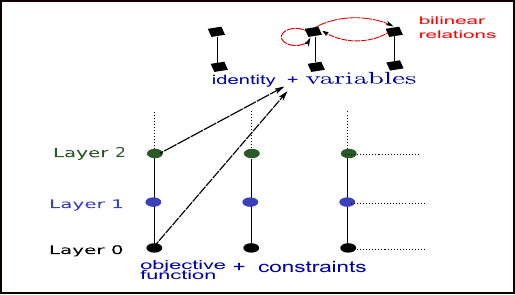}
    \caption{Illustration of \textbf{Graph 2}.}\label{fig:graph2}
    \end{figure}

%

\section{Computational Case}\label{sec:data_structures_computational_case}

The following example incorporates all the steps and the algorithms proposed in this paper. We construct the binary labelled graph and then enter it into \texttt{nauty} through dreadnaut command lines which compute the formulation group of the original problem.
\subsection{Numerical Example}
In this part we consider the Example \ref{ex1d} of a \ref{eq:QP} problem from Subsection \ref{sub:data_structures_matrix} and apply the two different graph representations as described in Subsection \ref{sub:graph}. 
Recall Example \ref{ex1d}:
\begin{align*} 
\max \quad & 3x_1 +  3x_4 + 2x_2x_3 \hspace{1.9cm} (c_0)\\
 	    & x_2 + {x_1}^2 + 1\leq 0\hspace{2.1cm} (c_1)\\
            & x_3 + x_4^2 + 1\leq 0 \hspace{2.27cm} (c_2)\\
            & x_2 + x_3 + 1\leq 0   \hspace{2.27cm} (c_3)\\
            & x_1, x_2, x_3, x_4 \in [0,1]         
\end{align*}
The sparse matrix representation:
\begin{itemize}
 \setlength{\itemsep}{0pt}
  \setlength{\parskip}{0pt}
  \setlength{\parsep}{0pt}
  \item $ \bf M_{QP_1} $ = $ (3 \, 3\, 2 \, 1\, 1\, 1\, 1\, 1\, 1\, 1\, 1\, 1)$
\item $     \bf I_{QP_1} $ = $ (0 \, 0 \, 2 \, 0  \, 0 \, 1 \, 0 \, 0 \, 4 \, 0 \, 0 \, 0 )$ 
\item $     \bf J_{QP_1} $ = $ (1  \,4 \, 3 \, 0 \, 2  \, 1  \, 0  \,  3 \,  4 \, 0 \, 2 \, 3)$
\item $     \bf K_{QP_1} $ = $ (0 \, 0 \, 0 \, 1\, 1\, 1\, 2\, 2\, 2\, 3\, 3\, 3)$
\end{itemize}
with vector of unique elements $\bf{U}$ $= (1\, 2\, 3)$. There are three unique elements and  $L = 2$ layers (see Equation \ref{eq:log}) to represent the relation of the variables of this problem. The binary representation of each unique element is computed using Equation \ref{eq:binary}, e.g.\ $3 = 2^1 + 2^0$ indicates that there is an edge between vertices on layer zero and another edge between the same vertices on layer 1. Following the Graph 2 description, this graph consists of 4 layers and $|\V| = 18$, one associated with a constant element and one with the objective function and the rest for the variables and constraints of the problem. The graph representation of \ref{ex1d} is shown in Figure 12. \texttt{Nauty} generates the following permutations: $\pi$ = (1 2)(5 6)(9 12)(10 11)(14 17)(15 16); the automorphism group of the graph under which it remains invariant. 
The relevant enumeration distinguishes which permutations are applied on the constraints and which on the variables of the problem. We then reflect these information on the original problem and explain its symmetric properties. Permutations (1 2)(5 6), permute the constraints  $c_1, c_2$ of Problem $QP_{1}$. Permutations (9 12)(10 11) are associated to the variables $x_1, x_4$ and $x_2, x_3$ with (14 17)(15, 16) their copies. Hence, the formulation group of problem \ref{ex1d} is  $\mathcal{G} = (x_1x_4)(x_2x_3)$. 

Problem $LPQ_1$ is the relaxed form of  the original Problem $QP_1$ after applying convex relaxation by introducing the auxiliary variables, $X_{23} = x_2x_3$, $X_{11} = x_1^2$, $X_{44} = x_4 ^2$ and adding the McCormick relaxation constraints.
\begin{align*} 
\max \quad & 3x_1 + 2X_{23} + 3x_4 \hspace{1.9cm} (c_0)\\
\text{s.t.} \quad & X_{11} + {x_2} + 1 \leq 0 \hspace{2cm} (c_1)\\
            & x_3 + X_{44} + 1\leq 0 \hspace{2cm} (c_2)\\
            & x_2 + x_3 + 1\leq 0 \hspace{2.2cm} (c_3)\\            		
	        & x_2 + x_3 - X_{23} - 1 \leq 0  \hspace{1.37cm} (c_4) \\ 
      		& X_{23} - x_2 \leq 0 \hspace{2.6cm} (c_5) \\
			& X_{23} - x_3 \leq 0 \hspace{2.6cm} (c_6) \\      		
      		& 2x_1 - X_{11} - 1\leq 0  \hspace{1.9cm} (c_7) \\ 
      		& X_{11} - x_1 \leq 0  \hspace{2.6cm} (c_8) \\      		
      		& 2x_4 - X_{44} - 1 \leq 0 \hspace{1.9cm} (c_9) \\ 
      		& X_{44} - x_4 \leq 0 \hspace{2.6cm} (c_{10}) \\
      		& X_{23}, X_{11}, X_{44} \geq 0 \quad \\ 
		& x_1, x_2, x_3, x_4 \in [0,1]  
\end{align*}
Similar to Problem $QP_1$ apply Method \ref{meth:1} described in Section \ref{sec:data_structures_BLG}. Observe that this problem also contains negative coefficients mapped to positive integers via function $\mu(z)$ e.g.\ $\mu(-1) = 4$ with binary representation $4 = 2^2$. The relevant graph is shown in Figure 13. \texttt{Nauty} generates (1 2)(5 6)(7 9)(8 10)(12 13)(16 17)(18 20)(19 21)(23 24)(27 28)(29 31)(30 32) (34 37)(35 36)(39 40) with specific permutations (34 37)(35 36)(39 40) to reveal the symmetric relations of variables $(x_1x_4)(x_2x_3)(X_{11}X_{44})$ the formulation group $\mathcal{G}$ of  $LPQ_1$. 
The above results validate both Method \ref{meth:1} and Method \ref{meth:2} in Section \ref{sec:data_structures_BLG} for representing an optimisation problem as a graph and then generate its symmetric properties. Figure \ref{fig:DAG_ex} shows a \textit{DAG} representation of Problem \ref{ex1d} as described in Section \ref{sec:data_structures_detect_symmetry_DAG}. The leaf nodes represent the variables and the coefficients, the intermediate nodes the operators and root nodes the plus signs that indicate the existence of a new constraint in this problem. Colours in the graph explain the vertex partitioning of the nodes that can be exchanged. The graph size in terms of the number of vertices and edges is smaller to the size of the methods proposed in this paper for this example and generates the same formulation group.
\begin{figure}
\centering
    \includegraphics[width=0.4\linewidth]{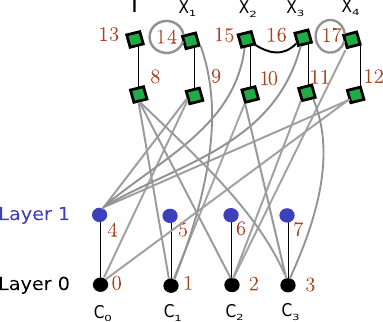}\label{fig:ex1}\caption{Illustration of Problem $QP_1$ using \textbf{Graph 2} representation.}
\end{figure}
\begin{figure}\centering
    \includegraphics[width=0.6\linewidth]{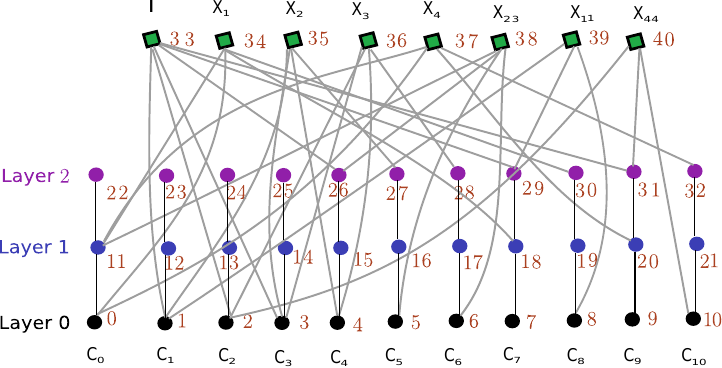}\label{fig:ex2}\caption{Illustration of Problem $LQP_1$ using \textbf{Graph 1} representation.}
\end{figure} 
\begin{figure}
\centering
    \includegraphics[width=0.6\linewidth]{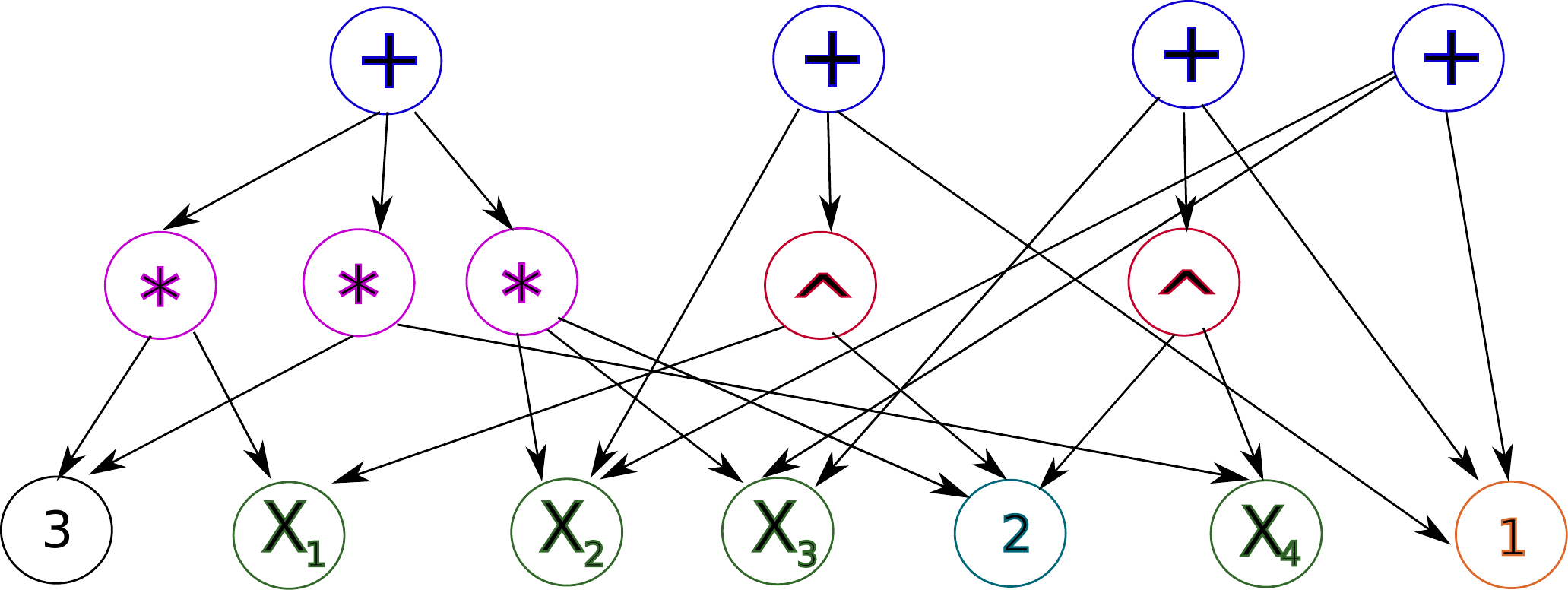}\caption{Illustration of Problem 1 using \textbf{DAG} representation.}\label{fig:DAG_ex}
\end{figure}

\subsection{Comparison with Current Methods}

We evaluate the trade-offs among the graph constructions in this paper and different graph constructions already in the literature. 
Regarding the graph transformation and its significant role in dealing with symmetry, \citet{ostrowski:2008} introduce the method "Orbital Branching" for combating symmetry. They illustrate each problem and its subproblems on each node of the tree as graphs and use \texttt{nauty} to compute the automorphism group and the orbits of the graph.
The presence of many coefficients in a problem expand the difficulty of identifying its symmetric properties.
 \citet{libertiref:2012} uses \textit{Directed Acyclic Graphs} to represent any mathematical expression of MINLP and automatically generates the formulation group.
 A major advance of both methods is that they are easy to implement and \textit{DAG} can capture the structure of any class of mathematical programming problem.
This work proposes an alternative method that may be useful when working with problems with many coefficients of different values.  
Using the function that assigns integer values to the coefficients of a problem let us work not only with non $0-1$ coefficient but with any other value.
Also the use of a logarithmic number of layers may reduce the number of nodes in the graph for problems with a large number of coefficients. 
For example, we are able to present 60 different coefficients in a graph with 6 layers. 
Another advantage is that we are able to capture the relation of bilinear terms in a way that it is unnecessary to create new nodes for every mathematical operation presented in the problem. Addition is the main operation on which the structure of the graph is based and multiplication is only presented with edges and loops. Subtraction is treated as a new coefficient together with the number that follows. 
The original form of \textit{DAG} graphs without any simplification, provides important informations on the exact formulation of the original problem something which is not clear with our methods.
\textit{BLG} may be associated with problems that have the exact same symmetric structure but different formulation. 
This work though focuses on providing an alternative method to detect the symmetric structure of a problem and not solving the problem itself.



\section{Conclusion}\label{sec:data_structures_conclusion}
This work appraise the presence and significance of symmetry in optimisation problems. Symmetry representation and detection are the fundamental steps towards exploiting symmetry. We propose graph structures that may capture the symmetric properties of a problem in a coherent size. 

\bibliographystyle{plainnat}
\bibliography{references}

\begin{thebibliography}{48}
\providecommand{\natexlab}[1]{#1}
\providecommand{\url}[1]{\texttt{#1}}
\expandafter\ifx\csname urlstyle\endcsname\relax
  \providecommand{\doi}[1]{doi: #1}\else
  \providecommand{\doi}{doi: \begingroup \urlstyle{rm}\Url}\fi

\bibitem[Alemany et~al.(2016)Alemany, Komarnicki, Linand, and
  Magnago]{magnago:2016}
J.~Alemany, P.~Komarnicki, J.~Linand, and F.~Magnago.
\newblock Exploiting symmetry in unit commitment solutions for a large-scale
  electricity market.
\newblock In \emph{Electric Power Systems Research}, volume 140. Elsevier,
  2016.

\bibitem[Anstreicher(2009)]{anstreicher:2009}
K.~M. Anstreicher.
\newblock Semidefinite programming versus the reformulation-linearization
  technique for nonconvex quadratically constrained quadratic programming.
\newblock \emph{J Glob Optim}, 43\penalty0 (2):\penalty0 471--484, 2009.

\bibitem[Armstrong(1988)]{armstrong:1988}
M.~A. Armstrong.
\newblock \emph{Groups and Symmetry}.
\newblock Undergraduate Texts in Mathematics. Springer, 1988.

\bibitem[Berthold and Pfetsch(2009)]{berthold:2008}
T.~Berthold and M.E. Pfetsch.
\newblock Detecting orbitopal symmetries.
\newblock In \emph{Proceedings of the Annual International Conference of the
  German Operations Research Society (GOR)}, pages 433--438. Springer Berlin
  Heidelberg, 2009.

\bibitem[B{\"o}di et~al.(2013)B{\"o}di, Herr, and Joswig]{bodi:2013}
R.~B{\"o}di, K.~Herr, and M.~Joswig.
\newblock Algorithms for highly symmetric linear and integer programs.
\newblock \emph{Math Program}, 137\penalty0 (1):\penalty0 65--90, 2013.

\bibitem[Bremner et~al.(2014)Bremner, Sikiri{\'c}, Pasechnik, Rehn, and
  Sch{\"u}rmann]{bremner:2014}
D.~Bremner, M.~S. Sikiri{\'c}, D.~V. Pasechnik, T.~Rehn, and A.~Sch{\"u}rmann.
\newblock Computing symmetry groups of polyhedra.
\newblock \emph{{LMS} J Comp Math}, 17\penalty0 (1):\penalty0 565--581, 2014.

\bibitem[Cameron(1999)]{cameron:1999}
P.~J. Cameron.
\newblock \emph{Permutation Groups}.
\newblock London Mathematical Society Student Texts. Cambridge uni., 1999.

\bibitem[Costa et~al.(2013)Costa, Hansen, and Liberti]{costa:2013}
A.~Costa, P.~Hansen, and L.~Liberti.
\newblock On the impact of symmetry-breaking constraints on spatial
  branch-and-bound for circle packing in a square.
\newblock \emph{Disc Appl Math}, 161\penalty0 (1):\penalty0 96--106, 2013.

\bibitem[Crawford et~al.(1996)Crawford, Ginsberg, Luks, and Roy]{Crawford:1996}
J.~M. Crawford, M.~L. Ginsberg, E.M. Luks, and A.~Roy.
\newblock Symmetry-breaking predicates for search problems.
\newblock In \emph{Proceedings of the Fifth International Conference on
  Principles of Knowledge Representation and Reasoning}, pages 148--159. Morgan
  Kaufmann Publishers Inc., 1996.

\bibitem[Dias and Liberti(2015)]{dias:2015}
G.~Dias and L.~Liberti.
\newblock Orbital independence in symmetric mathematical programs.
\newblock In \emph{Proceedings of the Combinatorial Optimization and
  Applications: 9th International Conference, {COCOA} 2015, TX, USA, 2015.},
  pages 467--480. Springer International Publishing, 2015.

\bibitem[Faenza and Kaibel(2009)]{kaibel:2009}
Y.~Faenza and V.~Kaibel.
\newblock Extended formulations for packing and partitioning orbitopes.
\newblock \emph{Math Oper Res}, 34\penalty0 (3):\penalty0 686--697, 2009.

\bibitem[Fischetti et~al.(2017)Fischetti, Liberti, Salvagnin, and
  Walsh]{fischetti:2017}
M.~Fischetti, L.~Liberti, D.~Salvagnin, and T.~Walsh.
\newblock Orbital shrinking: Theory and applications.
\newblock \emph{Disc. Appl. Math.}, 222:\penalty0 109--123, 2017.

\bibitem[Friedman(2007)]{friedman:2007}
E.~J. Friedman.
\newblock Fundamental domains for integer programs with symmetries.
\newblock In A.~W.~M. Dress, V.~Xu, and B.~Zhu, editors, \emph{{COCOA}},
  Lecture Notes in Computer Science, pages 146--153. Springer, 2007.

\bibitem[Gent and Smith(2000)]{gent:2000}
I.~P. Gent and B.~M. Smith.
\newblock Symmetry breaking in constraint programming.
\newblock In W.~Horn, editor, \emph{Proceedings of the {ECAI}-2000,14th
  European Conference on Artificial Intelligence, 2000}, pages 599--603. {IOS}
  Press, 2000.

\bibitem[Gent et~al.(2005)Gent, Kelsey, Linton, McDonald, Miguel, and
  Smith]{gent:2005}
I.~P. Gent, T.~Kelsey, S.~A. Linton, I.~McDonald, I.~Miguel, and B.~M. Smith.
\newblock Conditional symmetry breaking.
\newblock In P.~van Beek, editor, \emph{Proceedings of the Principles and
  Practice of Constraint Programming - CP 2005: 11th International Conference,
  Spain, 2005.} Springer Berlin Heidelberg, 2005.

\bibitem[Ghoniem and Sherali(2011)]{ghoniem:2011}
A.~Ghoniem and H.~D. Sherali.
\newblock Defeating symmetry in combinatorial optimization via objective
  perturbations and hierarchical constraints.
\newblock \emph{IIE Transactions}, 43\penalty0 (8):\penalty0 575--588, 2011.

\bibitem[Herr et~al.(2013)Herr, Rehn, and Sch{\"u}rmann]{herr:2013}
K.~Herr, T.~Rehn, and A.~Sch{\"u}rmann.
\newblock Exploiting symmetry in integer convex optimization using core points.
\newblock \emph{Oper Res Lett}, 41\penalty0 (3):\penalty0 298--304, 2013.

\bibitem[Hojny and Pfetsch(2015)]{hojny:2015}
C.~Hojny and M.~E. Pfetsch.
\newblock Symmetry handling via symmetry breaking polytopes.
\newblock In \emph{13th Cologne Twente Workshop on Graphs and Combinatorial
  Optimization, Istanbul, Turkey, 2015. proceedings}, pages 63--66, 2015.

\bibitem[Kaibel et~al.(2011)Kaibel, Peinhardt, and Pfetsch]{kaibel:2011}
V.~Kaibel, M.~Peinhardt, and M.~E. Pfetsch.
\newblock Orbitopal fixing.
\newblock \emph{Disc Optim}, 8\penalty0 (4):\penalty0 595--610, 2011.

\bibitem[Kallrath(2009)]{kallrath:2009}
J.~Kallrath.
\newblock Cutting circles and polygons from area-minimizing rectangles.
\newblock \emph{J Glob Optim}, 43:\penalty0 299--328, 2009.

\bibitem[Knueven et~al.(2017)Knueven, Ostrowski, and Pokutta]{knueven:2017}
B.~Knueven, J.~Ostrowski, and S.~Pokutta.
\newblock Detecting almost symmetries of graphs.
\newblock \emph{Math Program Comp.}, 2017.

\bibitem[Kouyialis and Misener(2017)]{kouyialis:2016}
G.~Kouyialis and R.~Misener.
\newblock Detecting symmetry in designing heat exchanger networks.
\newblock In \emph{Proceedings of the International Conference of Foundations
  of Computer-Aided Process Operations - {FOCAPO/CPC} 2017}, Arizona, 2017.

\bibitem[Kucherenko et~al.(2007)Kucherenko, Belotti, Liberti, and
  Maculan]{kucherenko:2007}
S.~Kucherenko, P.~Belotti, L.~Liberti, and N.~Maculan.
\newblock New formulations for the kissing number problem.
\newblock \emph{Disc. App. Math.}, 155\penalty0 (14):\penalty0 1837 -- 1841,
  2007.

\bibitem[Liberti(2004)]{Liberti:2004}
L.~Liberti.
\newblock Reformulation and convex relaxation techniques for global
  optimization.
\newblock \emph{Quarterly Journal of the Belgian, French and Italian Operations
  Research Societies}, 2\penalty0 (3):\penalty0 255--258, 2004.

\bibitem[Liberti(2008)]{liberti:2008}
L.~Liberti.
\newblock Automatic generation of symmetry-breaking constraints.
\newblock In B.~Yang, D.~Z. Du, and C.~A. Wang, editors, \emph{Combinatorial
  Optimization and Applications}, pages 328--338, Berlin, Heidelberg, 2008.
  Springer.

\bibitem[Liberti(2012{\natexlab{a}})]{liberti:2012}
L.~Liberti.
\newblock Symmetry in mathematical programming.
\newblock In J.~Lee and S.~Leyffer, editors, \emph{Mixed Integer Nonlinear
  Programming}, pages 263--283. Springer New York, 2012{\natexlab{a}}.

\bibitem[Liberti(2012{\natexlab{b}})]{libertiref:2012}
L.~Liberti.
\newblock Reformulations in mathematical programming: automatic symmetry
  detection and exploitation.
\newblock \emph{Math Program}, 131\penalty0 (1):\penalty0 273--304,
  2012{\natexlab{b}}.

\bibitem[Liberti and Ostrowski(2014)]{ostrowski:2014}
L.~Liberti and J.~Ostrowski.
\newblock Stabilizer-based symmetry breaking constraints for mathematical
  programs.
\newblock \emph{J Glob Optim}, 60\penalty0 (2):\penalty0 183--194, 2014.

\bibitem[Liberti et~al.(2010)Liberti, Cafieri, and Savourey]{rose}
L.~Liberti, S.~Cafieri, and D.~Savourey.
\newblock {T}he {R}eformulation-{O}ptimization {S}oftware {E}ngine.
\newblock In K.~Fukuda, J.~van~der Hoeven, M.~Joswig, and N.~Takayama, editors,
  \emph{Mathematical Software - {ICMS} 2010}, pages 303 -- 314. Springer Berlin
  Heidelberg, 2010.

\bibitem[Margot(2002)]{margot:2002}
F.~Margot.
\newblock Pruning by isomorphism in branch-and-cut.
\newblock \emph{Math Program}, 94\penalty0 (1):\penalty0 71--90, 2002.

\bibitem[Margot(2003{\natexlab{a}})]{margot:2003}
F.~Margot.
\newblock Small covering designs by branch-and-cut.
\newblock \emph{Math Program}, 94\penalty0 (2):\penalty0 207--220,
  2003{\natexlab{a}}.

\bibitem[Margot(2003{\natexlab{b}})]{margot:2003a}
F.~Margot.
\newblock Exploiting orbits in symmetric {ILP}.
\newblock \emph{Math Program}, 98\penalty0 (1):\penalty0 3--21,
  2003{\natexlab{b}}.

\bibitem[Margot(2007)]{margot:2007}
F.~Margot.
\newblock Symmetric {ILP}: Coloring and small integers.
\newblock \emph{Disc Optim}, 4\penalty0 (1):\penalty0 40--62, 2007.
\newblock Mixed Integer Programming.

\bibitem[Margot(2010)]{margot:2010}
F.~Margot.
\newblock Symmetry in integer linear programming.
\newblock In \emph{50 Years of Integer Programming 1958-2008: From the Early
  Years to the State-of-the-Art}, pages 647--686. Springer Berlin Heidelberg,
  2010.

\bibitem[McCormick(1976)]{mccormick:1976}
G.~P. McCormick.
\newblock Computability of global solutions to factorable nonconvex programs:
  Part i - convex underestimating problems.
\newblock \emph{Math Program}, 10\penalty0 (1):\penalty0 147--175, 1976.

\bibitem[McKay and Piperno(2014)]{nauty}
B.~D. McKay and A.~Piperno.
\newblock Practical graph isomorphism, ii.
\newblock \emph{J. Symb. Comp.}, 60:\penalty0 94--112, 2014.

\bibitem[Ostrowski et~al.(2008)Ostrowski, Linderoth, Rossi, and
  Smriglio]{ostrowski:2008}
J.~Ostrowski, J.~Linderoth, F.~Rossi, and S.~Smriglio.
\newblock Constraint orbital branching.
\newblock In \emph{Proceedings of the Integer Programming and Combinatorial
  Optimization: 13th International Conference, IPCO 2008 Bertinoro, Italy,
  2008.}, pages 225--239. Springer-Verlag, 2008.

\bibitem[Ostrowski et~al.(2010)Ostrowski, Vannelli, and Anjos]{anjos:2010}
J.~Ostrowski, A.~Vannelli, and M.~F. Anjos.
\newblock Symmetry in scheduling problems.
\newblock In \emph{Cahier du {GERAD} G-2010-69, GERAD,QC, Canada}, 2010.

\bibitem[Ostrowski et~al.(2011)Ostrowski, Linderoth, Rossi, and
  Smiriglio]{ostrowski:2009}
J.~Ostrowski, J.~Linderoth, F.~Rossi, and S.~Smiriglio.
\newblock Orbital branching.
\newblock \emph{Math Program}, 126\penalty0 (1):\penalty0 147--178, 2011.

\bibitem[Ostrowski et~al.(2015)Ostrowski, Anjos, and
  Vannelli]{ostrowski-etal:2015}
J.~Ostrowski, M.~F. Anjos, and A.~Vannelli.
\newblock Modified orbital branching for structured symmetry with an
  application to unit commitment.
\newblock \emph{Math Program}, 150\penalty0 (1):\penalty0 99--129, 2015.

\bibitem[Pfetsch and Rehn(2015)]{rehn:2015}
M.~E. Pfetsch and T.~Rehn.
\newblock A computational comparison of symmetry handling methods for mixed
  integer programs. technical report.
\newblock Technical report, Optimization Online, 2015.

\bibitem[Puget(2005)]{puget:2005}
J.~F. Puget.
\newblock Automatic detection of variable and value symmetries.
\newblock In P.~van Beek, editor, \emph{Proceedings of the Principles and
  Practice of Constraint Programming - {CP} 2005: 11th International
  Conference, Spain, 2005.}, pages 475--489. Springer Berlin Heidelberg, 2005.

\bibitem[Qualizza et~al.(2012)Qualizza, Belotti, and Margot]{qualizzaetal-2012}
A.~Qualizza, P.~Belotti, and F.~Margot.
\newblock Linear programming relaxations of quadratically constrained quadratic
  programs.
\newblock In L.~Lee and S.~Leyffer, editors, \emph{Mixed Integer Nonlinear
  Programming}, pages 407--426, New York, NY, 2012. Springer New York.

\bibitem[Ramani and Markov(2005)]{ramani:2005}
A.~Ramani and I.~L. Markov.
\newblock Automatically exploiting symmetries in constraint programming.
\newblock In B.~V. Faltings, A.~Petcu, F.~Fages, and F.~Rossi, editors,
  \emph{Recent Advances in Constraints}, pages 98--112. Springer-Verlag, 2005.

\bibitem[Ramani et~al.(2006)Ramani, Aloul, Markov, and Sakallah]{ramani:2006}
A.~Ramani, F.~Aloul, I.~Markov, and K.~A. Sakallah.
\newblock Breaking instance-independent symmetries in exact graph coloring.
\newblock \emph{JAIR}, 1:\penalty0 324-- 329, 03 2006.

\bibitem[Salvagnin(2005)]{salvagnin:2005}
D.~Salvagnin.
\newblock A dominance procedure for integer programming.
\newblock Master's thesis, University of Padova, 2005.

\bibitem[Sherali and Smith.(2001)]{sherali:2001}
H.~D. Sherali and J.~C. Smith.
\newblock Improving discrete model representations via symmetry considerations.
\newblock \emph{Manag. Sci.}, 47\penalty0 (10):\penalty0 1396--1407, 2001.

\bibitem[Specht(2018)]{pack}
E.~Specht.
\newblock \url{http://www.packomania.com/}, 2018.

\end{thebibliography}
\end{document}